\documentclass[sn-mathphys,Numbered]{sn-jnl}% Math and Physical 

\usepackage{graphicx}%
\usepackage{multirow}%
\usepackage{amsmath,amssymb,amsfonts}%
\usepackage{amsthm}%
\usepackage{mathrsfs}%
\usepackage[title]{appendix}%
\usepackage{xcolor}%
\usepackage{textcomp}%
\usepackage{manyfoot}%
\usepackage{booktabs}%
\usepackage{algorithm}%
\usepackage{algorithmicx}%
\usepackage{algpseudocode}%
\usepackage{listings}%

\usepackage{bm}
\usepackage[export]{adjustbox}
\usepackage{makecell}
\usepackage{caption, subcaption}
\usepackage[absolute,overlay]{textpos}
\usepackage[super]{nth}
\usepackage{array}
\newcolumntype{M}{ >{\centering\arraybackslash} m{7cm} }
\newcolumntype{C}{ >{\centering\arraybackslash} m{2.5cm} }

\newcommand{\lr}[1]{\left(#1\right)}
\newcommand{\abs}[1]{\left\lvert#1\right\rvert}
\newcommand{\norm}[1]{\left\lVert#1\right\rVert}
\newcommand{\uu}{\mathbf{u}}

\newcommand{\pp}{\mathbf{p}}
\newcommand{\qq}{\mathbf{q}}

\DeclareMathOperator*{\argmin}{arg\,min}

\begin{document}

% \title[Efficient Computation of Multiscale Hamiltonian Systems Aided by Machine Learning]{Efficient Computation of Multiscale Hamiltonian Systems Aided by Machine Learning}

\title{Stabilization of parareal algorithms for long time computation of a class of highly oscillatory Hamiltonian flows using data}

\author*[1]{\fnm{Rui} \sur{Fang}}\email{rfang@utexas.edu}

\author[2]{\fnm{Richard} \sur{Tsai}}\email{ytsai@math.utexas.edu}

\affil*[1]{\orgdiv{Oden Institute for Computational Engineering and Sciences}, \orgname{The University of Texas at Austin}, \orgaddress{\city{Austin}, \state{TX}, \postcode{78712}, \country{USA}}}

\affil[2]{\orgdiv{Department of Mathematics and Oden Institute for Computational Engineering and Sciences}, \orgname{The University of Texas at Austin}, \orgaddress{\city{Austin}, \state{TX}, \postcode{78712}, \country{USA}}}

\abstract{Applying parallel-in-time algorithms to multiscale Hamiltonian systems to obtain stable long time simulations is very challenging. In this paper, we present novel data-driven methods aimed at improving the standard parareal algorithm developed by Lion, Maday, and Turinici in 2001, for multiscale Hamiltonian systems. The first method involves constructing a correction operator to improve a given inaccurate coarse solver through solving a Procrustes problem using data collected online along parareal trajectories. The second method involves constructing an efficient, high-fidelity solver by a neural network trained with offline generated data. For the second method, we address the issues of effective data generation and proper loss function design based on the Hamiltonian function. We show proof-of-concept by applying the proposed methods to a Fermi-Pasta-Ulum (FPU) problem. The numerical results demonstrate that the Procrustes parareal method is able to produce solutions that are more stable in energy compared to the standard parareal. The neural network solver can achieve comparable or better runtime performance compared to numerical solvers of similar accuracy. When combined with the standard parareal algorithm, the improved neural network solutions are slightly more stable in energy than the improved numerical coarse solutions.}

\maketitle

\section{Introduction}

Hamiltonian systems are ubiquitous in astronomy, molecular dynamics, classical mechanics, and theoretical physics. We concentrate on the separable Hamiltonian case 
\begin{align}
    H\lr{\pp,\qq} = \frac{1}{2} \pp^T M^{-1} \pp + U(\qq), \quad \pp,\qq \in\mathbb{R}^{d} 
    \label{eq:separable_hamiltoniam}
\end{align}
where $\pp$ and $\qq$ are respectively the generalized momentum and position in a $d$-dimensional space, $M$ a diagonal matrix denoting the masses, and $U$ a smooth scalar function depending on the position $\qq$. Physically, the Hamiltonian is interpreted as the total energy of a system, consisting of the kinetic energy $K(\pp):= \frac{1}{2} \pp^T M^{-1} \pp$ and the potential energy $U(\qq)$. The dynamics of the system is given by Hamilton's equations
\begin{align}
    \dot{\qq}= H_\pp = M^{-1} \pp, \quad \dot{\pp} = - H_\qq = - \nabla U(\qq).
\end{align}

% An outstanding characteristic of Hamiltonian systems is the symplecticity of the flow. That is, the flow map $\Phi_{t}\lr{\qq,\pp}$ satisfies
% \begin{align}
%     \Phi_{t}^{\prime}\lr{\qq,\pp}^{T}J\Phi_{t}^{\prime}\lr{\qq,\pp}=J
% \end{align}
% for any fixed $t$. As a result, it is desired for numerical integration techniques of Hamiltonian systems to preserve the symplectic structure. 

Geometric integrators (methods that preserve geometric properties of the exact flow) such as the Velocity Verlet method are frequently used to simulate Hamiltonian flows \cite{geometric}.  It is proved that the preservation of geometric structures may greatly improve long-time numerical integration compared with general-purpose methods.

Even with the improved long-time stability and accuracy of geometric integrators, the computational complexity remains high for many physical applications. In particular, for systems with multiple time scales, accurate long-time integration is very difficult because the small stepsize required for stable integration of the fast motions lead to a large number of time steps. Computational multiscale algorithms aim at reducing computational complexity by exploiting the underlying multiscale structures. For example, for systems with sufficiently wide separation of scales and certain homogeneity and ergodic properties, the heterogeneous multiscale methods (HMM) can compute the effective systems with significantly reduced complexity \cite{engquist2005heterogeneous}. 
% \begin{align}
%     \frac{d}{dt} \mathbf{u} = f(\mathbf{u}; c_\epsilon)
% \end{align}

Recently, due to the increased availability of parallel processors in modern supercomputers, there has been rising interest in developing time-domain parallelization algorithms to reduce the wall-clock computation time for time-dependent problems. As the algorithm of choice, the work in this paper will rely on the parareal method, a parallel-in-time algorithm introduced by Lion, Turinici and Maday \cite{lionsparareal}. To set up, the time domain $[0,T]$ is divided into $N$ subintervals of length $\Delta t = T/N$. The parareal method involves two numerical solvers that advance a solution by $\Delta t$: an efficient but low-fidelity coarse solver, denoted by $C_{\Delta t}$, and an accurate but expensive fine solver, denoted by $F_{\Delta t}$. The coarse solver solutions are iteratively corrected by fine solver solutions computed on smaller time intervals in parallel. Formally, letting $\uu_n^{(k)}$ denote the solution computed at iteration $k$ and time $t_n = n\Delta t$, the parareal iterations are given by 
\begin{align}
    \uu^{(k+1)}_{n+1} &= C_{\Delta t} \uu^{(k+1)}_{n} + \lr{F_{\Delta t} \uu^{(k)}_{n} - C_{\Delta t} \uu^{(k)}_{n}}, \\
    \uu^{(k)}_{0} &=\uu_0, \quad \uu^{(0)}_{n+1} = C_{\Delta t} \uu^{(0)}_{n}, \quad k = 0, 1, 2, ...; n = 0, 1, ..., N-1. \nonumber
    \label{eq:parareal}
\end{align}
% See also \cite{ariel2016parareal} for a multiscale algorithm for oscillatory systems based on the parareal method. 
Ideally, with sufficiently accurate coarse solver, the iterations will quickly converge to the sequential solution computed by the fine solver $F_{\Delta t}^n \uu_0$. However, a closer analysis reveals that the convergence relies on the stability of the parareal iterations, and the standard parareal is only stable for dissipative problems \cite{gander2007analysis, ariel2016parareal}. For oscillatory and hyperbolic problems (such as Hamiltonian systems), the standard parareal scheme is known to perform badly because the convergence restricts the length of integration time \cite{gander2014analysis}. To be more specific, let $\mathbf{e}^{(k)}_{n} :=\uu^{(k)}_{n} - F_{\Delta t}^n \uu_0 $. The bound for the amplification $\mathbf{e}^{(k+1)}_{n}/\mathbf{e}^{(k)}_{n}$ depends on the sum $\sum_{j=0}^{n-k-2}{\norm{C_{\Delta t}}^j}$. For dissipative problems, $\norm{C_{\Delta t}} < 1$ and so the sum can be bounded by some constant independent of $n$. In contrast, for purely hyperbolic problems, $\norm{C_{\Delta t}}$ is close to 1 and the sum will be proportional to $n$, which causes the iterations to become unstable. 

Over the past years, many efforts have been devoted to developing stable parareal schemes for oscillatory and hyperbolic problems. For example, Dai et al. \cite{dai2013symmetric} proposed a symmetric variant of the parareal scheme and coupled it with projections to the constant energy manifold. % quote difficulties 
Another class of methods were developed based on the idea of making correction to parareal solutions using data. We will review more of these works in Section 2.

In this paper, we will focus on multiscale Hamiltonian systems where $H = H_\epsilon$. Here, $\epsilon$ is a small parameter indicating the small time or length scale.
% which arises when $U = U_0 +\epsilon^{-1} U_1$ or $M = M_\epsilon$, i.e., the particle masses vary by orders of magnitude and the smallest mass is at the order of $\epsilon$. % elliptic and hyperbolic U(q) 
Our aim is to improve the computational efficiency for long-time simulations of multiscale Hamiltonian systems by leveraging recent advancement in machine learning (ML) and parallel-in-time algorithms. Specifically, we would like to develop data-driven methods to stabilize the standard parareal scheme. Our idea is to use an improved solver $\Phi^{\theta, k}_{\Delta t}$ in place of $C_{\Delta t}$ in \eqref{eq:parareal}. The superscript $\theta$ in $\Phi^{\theta, k}_{\Delta t}$ denotes the unknown parameters defining the solver. The superscript $k$ indicates this solver may depend on the iteration. We propose two approaches to construct $\Phi^{\theta, k}_{\Delta t}$: 
\begin{enumerate}
    \item Enhance an existing coarse solver $C_{\Delta t}$ through a correction operator that aligns the ``phase'' in the coarse solutions $C_{\Delta t} \uu^{(k)}_n$ and fine solutions $F_{\Delta t} \uu^{(k)}_n$ for each iteration. In other words, $\Phi^{\theta, k}_{\Delta t} := \Psi^{(k)}_{\Delta t} \circ C_{\Delta t}  \approx F_{\Delta t}$, where $\Psi^{(k)}_{\Delta t}$ is the correction operator constructed from data collected \textit{online} during parareal iterations. Because $\Psi^{(k)}_{\Delta t}$ is obtained by solving an orthogonal Procrustes problem, this approach is named the Procrustes parareal method. The Procrustes parareal iteration is given by
    \begin{align}
        \uu^{(k+1)}_{n+1} = \Psi^{(k)}_{\Delta t} \circ C_{\Delta t}  \uu^{(k+1)}_{n} + \lr{F_{\Delta t} \uu^{(k)}_{n} - \Psi^{(k)}_{\Delta t} \circ C_{\Delta t} \uu^{(k)}_{n}}.
    \end{align}
        
    \item Approximate the fine solver (namely the solution map with a fixed stepsize $\Delta t$) using a neural network (NN), i.e., $\Phi^{\theta, k}_{\Delta t} := \Phi^{\text{NN}}_{\Delta t} \approx F_{\Delta t}$ where $\Phi^{\text{NN}}_{\Delta t}$ stands for a NN solver. Unlike in the Procruste parareal approach, the NN solver is constructed  using \textit{offline} training data. For this approach, we will address the issue of suitable training data generation and loss function design. The parareal iteration with an NN solver is as follows:  
    \begin{align}
        \uu^{(k+1)}_{n+1} = \Phi^{\text{NN}}_{\Delta t} \uu^{(k+1)}_{n} + \lr{F_{\Delta t} \uu^{(k)}_{n} -  \Phi^{\text{NN}}_{\Delta t} \uu^{(k)}_{n}}.
    \end{align}
    
    \end{enumerate}

% In principle, the two approaches can be easily combined, i.e., $\Phi^{\theta, k}_{\Delta t} := \Psi^{(k)}_{\Delta t} \circ \Phi^{\text{NN}}_{\Delta t} \approx F_{\Delta t}$. 

The paper is organized as follows. Section 2 presents the Procrustes parareal approach. Section 3 presents the neural network approach for approximating the solution map. Section 4 presents a case study for the Fermi-Pasta-Ulum (FPU) problem. We then conclude our findings in Section 5.

\section{Enhancing the coarse solvers by data}

% \footnotetext[1]{Similar works: \cite{yalla2018parallel}}

In order to address the instability issue in the standard parareal iterations, several methods involving a correction to bridge the gap between fine and coarse solutions were proposed in the past. For example, Farhat and Chandesris \cite{farhat2003time} used a Newton-type iteration to reduce the jumps between the fine and coarse solutions. In \cite{ariel2016parareal}, Ariel et al. proposed the $\theta$-parareal scheme, which uses an interpolation-based linear operator to enhance the coarse solver for oscillatory systems. In \cite{nguyen2020stable}, Nguyen and Tsai focused on the second-order wave equations and developed a correction operator based on minimizing the wave energy residual of the fine and coarse solutions. The resulting correction successfully stabilizes the parareal iterations by aligning the ``phase'' in the wave fields computed by the fine and coarse solver respectively. Later, the same authors proposed in \cite{nguyen2023numerical} a deep learning approach to enhance the coarse solver to reduce the ``phase'' errors in wave propagation. 

The success for the wave equations in \cite{nguyen2020stable} directly inspires the development of a similar approach for a class of Hamiltonian systems where the notion of ``phase'' can be suitably defined. 

In this section, we first provide our definition of ``phase'' for Hamiltonian systems, and then introduce the procedure to obtain a correction operator $\Psi^{(k)}_{\Delta t}$ from data computed in previous iterations $$\{ C_{\Delta t} \uu^{(k^\prime)}_n, F_{\Delta t} \uu^{(k^\prime)}_n\} \quad  n=0, 1, \cdots, N-1; \; k^\prime=0, 1, \cdots, k-1.$$

\subsection{A practical notion of ``phase'' for Hamiltonian systems}

For integrable Hamiltonian systems such as harmonic oscillators or Kepler systems, the ``phase'' can be naturally defined as the angle variables from the action-angle coordinates. For non-integrable systems, such as the FPU problem,  where action-angle coordinates are not available, we need alternative definitions for ``phase''.  

Because phase is an angle-like object, for a separable Hamiltonian~\eqref{eq:separable_hamiltoniam}, it is natural to consider a transform function $\Lambda$, which maps the energy level set to a hypersphere, whose radius satisfies 
\begin{align}
    \norm{\Lambda\lr{[\pp, \qq]}}_2^2=H\lr{\pp, \qq} + \text{constant}.
\end{align}
This way, we can define the ``phase'' difference between $[\pp_1, \qq_1]$ and $[\pp_2, \qq_2]$ as the angle between the transformed vectors $\Lambda\lr{[\pp_1, \qq_1]}$ and $\Lambda\lr{[\pp_2, \qq_2]}$. We call $\Lambda$ the energy transform because the $l_2$ norm of the transformed vector is related to the energy. 

% therefore it makes sense to compare using Procrustes 

% We define the energy semi-norm 
% \begin{align}
%     E \lr{\begin{bmatrix}\pp \\ \qq\end{bmatrix}} := \norm{\Lambda \lr{ \begin{bmatrix}\pp \\ \qq\end{bmatrix}}}_2^2 = H(\pp, \qq).
%     \label{eq:energyseminorm}
% \end{align}
% Here, the transform function $\Lambda$ maps the phase space variables $(\pp,\qq)$ to energy components $\mathbf{z}$
The specific form of $\Lambda$ and the pseudo-inverse $\Lambda^\dagger$ depend on the Hamiltonian. Noticeably, not all Hamiltonian functions allow a valid definition of $\Lambda$. 
% In particular, $\Lambda$ is only defined when the level set of $H\lr{\pp, \qq}$ is diffeomorphic to a sphere. 
For example, for a 1D harmonic oscillator whose Hamiltonian is $H(p,q) = \frac{1}{2}p^2 + \frac{1}{2}q^2$, 
\begin{align}
    \Lambda\lr{\begin{bmatrix}p \\ q\end{bmatrix}} := \begin{bmatrix}p / \sqrt{2} \\ q /\sqrt{2} \end{bmatrix}, \quad 
    \Lambda^{\dagger}\lr{\begin{bmatrix}\Tilde{p} \\ \Tilde{q}\end{bmatrix}} := \begin{bmatrix}\sqrt{2} \Tilde{p}\\ \sqrt{2} \Tilde{q}\end{bmatrix}.
\end{align}
For a 2D Kepler problem where $H(\pp, \qq) = \frac{1}{2} \pp^T \pp - \frac{1}{\norm{\qq}}$, $\Lambda$ does not exist because of the negative potential term.

In this paper, we will work with the FPU problem and will show the corresponding definition of $\Lambda$ in Section 4.

\subsection{The Procrustes Parareal}

% the wave energy semi-norm (i.e. $\Lambda$ transforms the wave field solution $\uu = (u, u_t)$ to the energy components $(\nabla u, c^{-1}      u_t)$ and $\norm{}_* = \norm{}_2$) and seek to find a phase corrector by minimizing the wave energy residual of the fine and coarse solutions. The correction operator is defined as $\Psi^{(k)}_{\Delta t} := \Lambda^{\dagger} \Omega^{(k)}_{\Delta t} \Lambda$ where $\Lambda^{\dagger}$ is the pseudo-inverse of $\Lambda$, and $\Omega^{(k)}_{\Delta t}$ is a unitary matrix solving \eqref{eq:min_problem_online_data} (thus equivalent to a Procrustes problem). 

% For a general problem, the correction operator will minimize the following difference 
% \begin{align}
%     \Psi^{(k)}_{\Delta t} = \argmin_{ \Psi_{\Delta t} \in X} \sum_{k^{\prime}=0}^{k} \alpha^{k^{\prime}} \sum_{n=0}^{N-1}  \norm{ \Lambda \lr{F_{\Delta t} \uu^{(k^{\prime})}_n} - \Lambda \lr{ \Psi_{\Delta t} C_{\Delta t} \uu^{(k^{\prime})}_n}}_*
%     \label{eq:min_problem_online_data}
% \end{align}
% where $X$ is a functional space in which one searches for $\Psi_{\Delta t}$, $\alpha^{k^{\prime}}\in(0, 1]$ is the weight put to different iterations, $\Lambda$ is a suitable transform applied to the data before taking differences and $\norm{}_*$ is a suitable norm for measuring the differences. 

In the following, we use $\uu$ to represent the concatenated vector $[\pp, \qq]$.  Assuming $\Lambda$ and $\Lambda^{\dagger}$ are both known, we define the correction operator as
\begin{align}
    \Psi^{(k)}_{\Delta t} := \Lambda^{\dagger} \circ \Omega^{(k)}_{\Delta t} \circ \Lambda,
\end{align}
where $\Omega^{(k)}_{\Delta t}$ is an orthogonal transformation to be determined. We see the advantage of defining $\Lambda(\uu)$ - any orthogonal transformation on $\Lambda(\uu)$ preserves $H(\uu)$. Hence, the correction operator preserves the energy. 

The Procrustes parareal method is given as follows: (function composition symbols are left out for brevity) 
\begin{align}
    \uu^{(k+1)}_{n+1} &= \Lambda^{\dagger} \Omega^{(k)}_{\Delta t}  \Lambda C_{\Delta t} \uu^{(k+1)}_{n} + \lr{F_{\Delta t} \uu^{(k)}_{n} - \Lambda^{\dagger} \Omega^{(k)}_{\Delta t}  \Lambda C_{\Delta t} \uu^{(k)}_{n}},  \\
    \uu^{(k)}_{0} &=\uu_0, \quad \uu^{(0)}_{n+1} = C_{\Delta t} \uu^{(0)}_{n}, \quad k = 0, 1, 2, ...; n = 0, 1, ..., N-1. \nonumber
\end{align}
Here, $\Omega^{(k)}_{\Delta t}$ is obtained by solving the orthogonal Procrustes problem 
\begin{align}
    \Omega^{(k)}_{\Delta t} &:= \argmin_\Omega \sum_{n=0}^{N-1}  \norm{ f_n - \Omega g_n}_2^2 \quad \text{s.t.} \quad \Omega\Omega^T = \Omega^T\Omega = I, \\
    f_n &:= \Lambda F_{\Delta t} \uu^{(k)}_n, \nonumber \\
    g_n &:= \Lambda C_{\Delta t} \uu^{(k)}_n. \nonumber
\end{align}

The geometric interpretation is to minimize the sum of the phase errors between fine and coarse solutions computed along the trajectory from the last iteration. Hence $\Omega^{(k)}_{\Delta t}$ is referred to as the phase corrector. 

% Let $\alpha^{k^{\prime}} = 
%     \begin{cases}
%     1 & k^\prime = k \\
%     0 & k^\prime < k
%     \end{cases}$.    % this is the current version. should probably try \alpha^k != 0 for k < 0 

We follow the standard way to solve the orthogonal Procrustes problem, which uses the singular value decomposition (SVD) of the correlation matrix $M:=F G^T$. Here, $F := [f_0 \; f_1 \; \cdots \; f_{N-1}]$, $G := [g_0 \; g_1 \; \cdots \; g_{N-1}]$. Let $U\Sigma V^T$  be the SVD of $M$. If $M$ has full rank, then the minimizer is uniquely $\Omega_* = UV^T$. We refer readers to \cite{gower2004procrustes} for more details of the Procrustes problem. 

The pseudo-code for the Procrustes parareal method is provided in Algorithm \ref{alg:procrustes_parareal}.

\begin{algorithm}
\caption{Procrustes parareal method}\label{alg:procrustes_parareal}
\begin{algorithmic}
\Require{initial state $\uu_0$, number of time intervals $N$, number of iterations $K$, coarse solver $C_{\Delta t}$, fine solver $F_{\Delta t}$, map from phase space variables to energy components $\Lambda$, map from energy components to phase space variables $\Lambda^{\dagger}$}
\Ensure{solutions $\uu^{(k)}_{n}$ for $k=0, 1, ...; n=0, 1, ..., N$}
\State Compute initial solution:
\State $\uu^{(0)}_{0} =\uu_0$
\For{$n = 0:N-1$}
    \State $\uu^{(0)}_{n+1} = C_{\Delta t} \uu^{(0)}_{n}$
\EndFor
\State 
\State $k=0$
\While{$k \leq K$}
    \State Compute coarse and fine solutions in parallel:
    \For{$n = 0:N-1$} 
        \State $f_n = F_{\Delta t} \uu^{(k)}_n$
        \State $g_n = C_{\Delta t} \uu^{(k)}_n$ 
    \EndFor
    \State 
    \State Construct data matrices and solve the Procrustes problem:
    \State $F = [\Lambda(f_0)\; \Lambda(f_1)\; \cdots \; \Lambda(f_{N-1})]$  
    \State $G = [\Lambda(g_0)\; \Lambda(g_1)\; \cdots \; \Lambda(g_{N-1})]$
    \State $M = F G^T$
    \State $[U, \Sigma, V] = SVD(M)$
    \State $\Omega = U V^T$
    \State
    \State Compute solution of current iteration sequentially:
    \State $\uu^{(k+1)}_{0} =\uu_0$ 
    \For{$n = 0:N-1$}
       \State $\uu^{(k+1)}_{n+1} = \Lambda^{\dagger} \Omega  \Lambda C_{\Delta t} \uu^{(k+1)}_{n} + \lr{F_{\Delta t} \uu^{(k)}_{n} - \Lambda^{\dagger} \Omega \Lambda C_{\Delta t} \uu^{(k)}_{n}}$
    \EndFor
    \State $k=k+1$
\EndWhile
\end{algorithmic}
\end{algorithm}

%[properties of Procrustes parareal]
% potential weaknesses of procrustes parareal 
% relies on span of previous iteration solutions 
% corrected solution 

\section{Neural network approximation of the solution map}

In this section, we present our second data-driven approach, which is to use a neural network (NN) to approximate the solution map $F_{\Delta t}$. To construct the NN solver $\Phi^{\text{NN}}_{\Delta t}$, the main task is to solve an optimization problem
\begin{align}
    \Phi^{\text{NN}}_{\Delta t} = \argmin_{\Phi^{\text{NN}}_{\Delta t} \in X} \frac{1}{\lvert \mathcal{D}_0 \rvert} \sum_{\uu_0 \in \mathcal{D}_0}  l\lr{\uu_0, \Phi^{\text{NN}}_{\Delta t}, F_{\Delta t}},
    \label{eq:min_problem_offline_data}
\end{align}
where $X$ is a function space determined by the network architecture, $\mathcal{D}_0$ a set of input data points, and $l\lr{\uu_0, \Phi^{\text{NN}}_{\Delta t}, F_{\Delta t}}$ the misfit term for each data point $\uu_0$ given the reference solution map $F_{\Delta t}$ and the approximated map $\Phi^{\text{NN}}_{\Delta t}$. We describe our setup for each of these components as follows.

\subsection{Choice of neural network architecture}

For simplicity, we choose fully-connected residual networks (ResNets). Compared to a regular multilayer perceptron, the residual network adds skip connections between pairs of hidden layers. 
%motivation? 

Let $L$ denote the number of hidden layers and $n$ the number of nodes per hidden layer. The layer outputs are defined as follows:
\begin{align}
    \text{input layer:} \quad y^{(0)} &:= x \in \mathbb{R}^{2d}, \nonumber \\
    \text{1st hidden layer:} \quad y^{(1)} &:= \sigma(W^{(1)}y^{(0)} + b^{(1)}) \in \mathbb{R}^{n}, \nonumber \\
    \text{$l$-th hidden layer:} \quad y^{(l)} &:= y^{(l-1)} + \frac{1}{L} \sigma(W^{(l)}y^{(l-1)} + b^{(l)}) \in \mathbb{R}^{n}, \quad l=2,...,L \nonumber \\
    \text{output layer:} \quad y^{(L+1)} &:= W^{(L+1)}y^{(L)} + b^{(L+1)} \in \mathbb{R}^{2d}.
\end{align}
Here, $W^{(l)}$ and $b^{(l)}, l=1,...,L+1$ are weights and biases to be determined through the training procedure. We use the Exponential Linear Unit (ELU) \cite{clevert2015fast} for the nonlinear activation function $\sigma$. Note that we adopt a scaling factor $1/L$ for the hidden layers with skip connections. This technique was proposed in \cite{weinan2020machine} to make the network performance more robust against hyperparameter change.  

% Hamiltonian Reversible Network (HRN) \cite{chang2018reversible}   
%     \begin{align*}
%         y^{(l+1)} &= 
%         \begin{bmatrix}
%             p^{(l+1)} \\
%             q^{(l+1)}
%         \end{bmatrix} := 
%         \begin{bmatrix}
%             p^{(l)} + h (K_1^{(l+1)})^{T} \sigma(K_1^{(l+1)}q^{(l)} + b_1^{(l+1)})\\
%             q^{(l)} - h (K_2^{(l+1)})^{T} \sigma(K_2^{(l+1)}p^{(l+1)} + b_2^{(l+1)})
%         \end{bmatrix}
%     \end{align*}
% Advantage: reversible

\subsection{Design of misfit term}

The function to be approximated is a solution map that maps a phase space state to another state. This allows us to use a sequence of successive time steps to construct the misfit term. Suppose that for an input $\uu_0$ and a sequence length $S$, we generate $\{ \uu_i \}_{1\leq i \leq S}, \uu_i = \lr{F_{\Delta t}}^{i} \uu_0 $ as the target sequence and $\{ \tilde{\uu}_i \}_{1\leq i \leq S}, \tilde{\uu}_i = \lr{\Phi^{\text{NN}}_{\Delta t}}^{i} \uu_0 $ as the approximated sequence. The misfit is then computed between the approximated sequence and the target sequence
\begin{align}
    l\lr{\uu_0, \Phi^{\text{NN}}_{\Delta t}, F_{\Delta t}} = \frac{1}{S} \sum_{i=1}^{S} \text{diff} \lr{\uu_i, \tilde{\uu}_i}.
\end{align}
We remark that because the network is applied recursively for obtaining $\tilde{\uu}_i$, this multi-step loss essentially makes training the network like training a recurrent neural network.

There are several ways to measure the difference between $(\uu_i, \tilde{\uu}_i)$. One common approach is to use the Euclidean metric of $\mathbb{R}^{2d}$. This is known as the mean squared error 
\begin{align}
    \text{MSE}\lr{\uu_i, \tilde{\uu}_i} = \norm{\uu_i -\tilde{\uu}_i}_2^2.
\end{align} 
The Euclidean metric puts equal weights on $\pp$ and $\qq$ components of $\uu$. While minimizing the mean squared error aligns with the goal of reducing the trajectory error, it often leads to imbalanced energy error because the Hamiltonian function does not always weight $\pp$ and $\qq$ similarly. Therefore, to naturally balance the components based on the Hamiltonian, we adopt the energy transform $\Lambda$ as defined in Section 2.1. We define the energy balanced error 
\begin{align}
    \text{EBE}\lr{\uu_i, \tilde{\uu}_i} = \norm{\Lambda \uu_i - \Lambda \tilde{\uu}_i}_2^2.
\end{align}

To put together, suppose we use the energy balanced error, the misfit term for an initial state $\uu_0$ and a sequence length $S$ is given by
\begin{align}
    l\lr{\uu_0, \Phi^{\text{NN}}_{\Delta t}, F_{\Delta t}} = \frac{1}{S} \sum_{i=1}^{S} \norm{\Lambda \lr{\lr{F_{\Delta t}}^{i} \uu_0} - \Lambda \lr{\lr{\Phi^{\text{NN}}_{\Delta t}}^i \uu_0}}_2^2.
\end{align}

\subsection{Generation of input data set}

The next problem is to generate a proper set of initial conditions $\uu_0$ for training. Unlike in the Procrustes parareal approach, where the data are collected online during parareal iterations, here we have to generate training data offline to construct an effective NN solver.  

In order to understand ``what is an reasonable distribution to sample in the phase space'', we regard the misfit term in the optimization problem \eqref{eq:min_problem_offline_data} as the mean error over a continuous distribution of $\uu$ in the phase space:
\begin{align}
    \frac{1}{\lvert \mathcal{D}_0 \rvert} \sum_{\uu \in \mathcal{D}_0}  l\lr{\uu, \Phi^{\text{NN}}_{\Delta t}, F_{\Delta t}}  \approx \int_{\mathbb{R}^{2d}} l\lr{\uu, \Phi^{\text{NN}}_{\Delta t}, F_{\Delta t}}  d\mu(\uu). 
\end{align}
It is thus natural to consider using a relevant invariant measure of the Hamiltonian flow for $\mu$. 

Suppose we are interested in simulations of the Hamiltonian flow with a fixed total energy. Then, we should consider sampling 
%Since the total energy is conserved under a Hamiltonian flow, we 
an invariant measure on an energy level set
%, i.e., $D = \mathcal{M}_{H_0}$ where 
\begin{align}
    \mathcal{M}_{H_0} := \left \{ (\pp,\qq) \in\mathbb{R}^{2d} \mid H(\pp,\qq) = H_0\right \}.
\end{align}
% Before introducing our sampling algorithms, we first review some numerical considerations. Although the true solution map preserves $H$ exactly, 
However, the numerical approximations may not preserve the total energy.
To start, the accurate fine solver used to approximate the true solution map is a symplectic integrator for which exact energy preservation is not possible. In addition, there is no guarantee for energy preservation by the general NN solver considered in this paper. 

Notice that by Liouville's theorem, the Hamiltonian flow preserves phase space volume. Then, we can construct an invariant measure in $\mathbb{R}^{2d}$ that concentrates on the chosen energy level $H_0$ as follows, using the coarea formula:
\begin{equation}
\begin{aligned}
&\int_{\mathbb{R}^{2d}} l\lr{\uu, \Phi^{\text{NN}}_{\Delta t}, F_{\Delta t}}  \exp^{-(H(\uu)-H_0)^2/2\sigma^2}d\uu\\
&= \int_{\mathbb{R}} \exp^{-(E-H_0)^2/2\sigma^2} \int_{\{H(\uu)=E\}} l\lr{\uu, \Phi^{\text{NN}}_{\Delta t}, F_{\Delta t}} \frac{dS}{|\nabla H(\uu)|}  dE. 
\end{aligned}
\end{equation}
In other words, one can separately sample the invariant densities on the energy level sets and the Gaussian density in the normal directions on the energy level sets.
% Hence, we propose that the training data to be sampled from a neighborhood of $\mathcal{M}_{H_0}$:
% \begin{align}
%     \mathcal{M}^\sigma_{H_0} := \bigcup_i \mathcal{M}_{H^i} \quad \text{where} \; H^i \sim \mathcal{N}(H_0, \sigma^2).
% \end{align}
% To sample from $\mathcal{M}^\sigma_{H_0}$, 

Motivated by this observation, we propose a novel sampling algorithm called HMC-$H_0$. The name comes from its resemblance to the Hamiltonian Monte-Carlo (HMC) algorithm~\cite{DUANE1987216}. Starting with $\qq=\qq_0$, we generate a chain of points by repeating the following two steps: 
\begin{enumerate}
    \item Momentum refreshment: randomly sample $\pp$ from the hypersphere defined by
    \begin{align}
        \left\{ \pp\in\mathbb{R}^d \mid \pp^T M^{-1} \pp = 2 \lr{ H_0 - U(\qq) }\right\};
    \end{align}
    \item Time integration: $\lr{\pp, \qq } \leftarrow F_{\delta t}\lr{\pp, \qq }$.
\end{enumerate}
Note that our approach is different from the original HMC algorithm in the momentum refreshment step. In the original HMC, $\pp$ is randomly sampled from a Gaussian distribution independent of current $\qq$, whereas in our approach, the distribution depends on $\qq$ and a fixed $H_0$.

% Integration time $\delta t$ must be carefully chosen
% (what is the resulting distribution?) 
% Let $\rho(\pp,\qq)$ be the density of $\mu_{\mathcal{M}_{H_0}}$. We can write $\rho(\pp,\qq) = \rho(\qq \mid \pp)\rho(\pp)$.  

%Now we briefly return to the discussion on invariant measures on an energy level set.
% One example of a flow-invariant measure is the Liouville measure. 

%The Liouville density is constant on an energy level set, according to the Liouville theorem. 
%Therefore, it is reasonable to that the training data should be sampled from this measure.
%This is a nontrivial task mainly due to the form of the nonlinear potential energy $U(\mathbf{q})$. 
We shall compare the HMC-$H_0$ algorithm to the following naive approach, combining random sampling in the momentum space and generating trajectories in the phase space using the flow. 
We first sample a set of momenta, then using the sampled momenta and $\qq_0$, we generate an ensemble of trajectories by flowing the points for a duration of time. The points along the trajectories are collected. This is a naïve attempt to sample the Liouville measures on the energy level sets, assuming ergodicity of the flows. We call this algorithm TrajEnsemble-$H_0$.

% If $L=1$, in the limit $\delta t \to 0, n \to \infty, m \to \infty$, the sampled distribution approaches $\mu_{\mathcal{M}^{\sigma}_{H_0}}$ [proof?]

Full descriptions of the algorithms are given in Algorithm~\ref{alg:hmc-H0} and Algorithm~\ref{alg:trajensemble-H0}. For both algorithms, we can leverage parallel computation to obtain a large number of data samples.

\begin{algorithm}
\caption{HMC-$H_0$}\label{alg:hmc-H0}
\begin{algorithmic}
\Require{target energy $H_0$, initial position $\qq_0 \in \mathbb{R}^d$, standard deviation $\sigma$, mass matrix $M$, number of chains $N_{\mathrm{chains}}$, number of transitions per chain $N_{\mathrm{trans}}$, integration time $\delta t$, reference solution map $F_{\delta t}$}
\Ensure{a set of phase space points $\{ (\pp^{i,j}, \qq^{i,j}) \}_{1\leq i \leq N_{\mathrm{chains}}, 1\leq j \leq N_{\mathrm{trans}}}$}
\For{$i = 1:N_{\mathrm{chains}}$}
    \State $\qq^{i,0} = \qq_0$
    \For{$j = 1:N_{\mathrm{trans}}$}
        \State Step 1: Momentum refreshment
        \State $K^{\prime} = 0$
        \While{$K^{\prime} \leq 0$}
            \State sample $K^{\prime} \sim \mathcal{N}(H_0-U(\qq^{i,j-1}), \sigma^2)$
        \EndWhile
        \State sample $\Tilde{\pp}^{\prime} \sim \mathcal{U}(\mathbb{S}^{d-1})$
        \State $\pp^{\prime} = \sqrt{2 K^{\prime}} M^{1/2}\Tilde{\pp}^{\prime}$
        \State 
        \State Step 2: Time integration
        \State $\pp^{i,j}, \qq^{i,j} = F_{\delta t}(\pp^{\prime}, \qq^{i,j-1})$
    \EndFor
\EndFor
\end{algorithmic}
\end{algorithm}

\begin{algorithm}
\caption{TrajEnsemble-$H_0$}\label{alg:trajensemble-H0}
\begin{algorithmic}
    
\Require{target energy $H_0$, initial position $\qq_0 \in \mathbb{R}^d$, standard deviation $\sigma$, mass matrix $M$, number of energy level sets $N_{\mathrm{level sets}}$, number of trajectories per level set $N_{\mathrm{traj}}$, length of trajectory $L$, step size of trajectory $\delta t$, reference solution map $F_{\delta t}$}
\Ensure{a set of phase space points $\{ (\pp^{i,j,k}, \qq^{i,j,k}) \}_{1\leq i \leq N_{\mathrm{level sets}}, 1\leq j \leq N_{\mathrm{traj}}, 1\leq k \leq L}$}
\For{$i = 1:N_{\mathrm{level sets}}$}
    \State $K^{i} = 0$
    \While{$K^{i} \leq 0$}
        \State sample $K^{i} \sim \mathcal{N}(H_0-U(\qq_0), \sigma^2)$
    \EndWhile
    \For{$j = 1:N_{\mathrm{traj}}$}
        \State sample $\Tilde{\pp}^{i,j} \sim \mathcal{U}(\mathbb{S}^{d-1})$
        \State $\pp^{i,j} = \sqrt{2 K^i} M^{1/2}\Tilde{\pp}^{i,j}$
        \For{$k=1:L$}
            \State $\pp^{i,j, k}, \qq^{i,j, k} = F_{\delta t}^k (\pp^{i,j}, \qq_0)$
        \EndFor
    \EndFor
\EndFor
\end{algorithmic}
\end{algorithm}

\section{Case study: The Fermi-Pasta-Ulum Problem}

We consider the Fermi-Pasta-Ulam (FPU) problem as a model problem to demonstrate properties of our proposed methods. 

First studied in 1955, the FPU problem \cite{fermi1955studies} describes a simple yet important model for nonlinear physics, which exhibits unexpected dynamical behaviors after long enough integration time. The model involves a chain of particles connected by springs that obey Hooke's law but with a weak nonlinear perturbation. Here, we adopt a version of the problem presented in \cite{geometric}. Suppose there are $2m$ mass points connected by alternating soft nonlinear springs and stiff linear springs. The variables $q_1, \cdots, q_{2m}$ ($q_0 = q_{2m+1} = 0$) represent the displacements of the mass points from equilibrium, and $p_i = d q_i/dt$ represent velocities. The Hamiltonian of this system is given by
\begin{align}
    H(\pp,\qq) = \frac{1}{2} \sum_{i=1}^{m} \lr{p_{2i-1}^2 + p_{2i}^{2}} + \frac{\omega^2}{4} \sum_{i=1}^{m} \lr{q_{2i}-q_{2i-1}}^2 +\sum_{i=0}^{m} \lr{q_{2i+1}-q_{2i}}^4,
    \label{eq:H_fpu}
\end{align}
where $\omega \gg 1$ is the frequency of the stiff linear springs. 

The dynamics of such a system has different behaviors on several different time scales. On the smallest time scale $\mathcal{O}(\omega^{-1})$, the linear springs show almost-harmonic oscillations with period close to $\pi/\omega$. On the time scale $\mathcal{O}(\omega^{0})$, the motion of the nonlinear springs becomes apparent. On the time scale $\mathcal{O}(\omega)$, there is slow energy exchange among the stiff springs. An illustration of motion on different time scales can be found in Section XIII.2 in \cite{geometric}. 

In our experiments, we aim to obtain stable simulation of the system on a time scale of $\mathcal{O}(\omega)$ using solvers with $\Delta t=1.0$. We will use $m=3$ (hence the degree of freedom $d=6$) and $\omega=300$. This is a challenging regime because the separation of characteristic time scales is large, i.e. from $1/300$ to 300. We will run the algorithms from the initial condition 
\begin{align}
    \pp_{\text{init}} = \begin{bmatrix}
        0 & \sqrt{2} & 0 & 0 & 0 & 0
    \end{bmatrix}^T, \quad 
    \qq_{\text{init}} = \begin{bmatrix}
        \frac{1-\omega^{-1}}{\sqrt{2}} & \frac{1+\omega^{-1}}{\sqrt{2}} & 0 & 0 & 0 & 0
    \end{bmatrix}^T.
    \label{eq:test_u0}
\end{align}
The corresponding energy is 
\begin{align}
    H(\pp_{\text{init}},\qq_{\text{init}}) = 2 + 3\omega^{-2} + \frac{1}{2}\omega^{-4}.
    \label{eq:test_H0}
\end{align}

Given a reference solution $\uu^{\text{ref}}=\lr{ \pp^{\text{ref}}, \qq^{\text{ref}}}$ and a computed solution $\uu=\lr{\pp,\qq}$, we shall report the trajectory error 
\begin{align}
    \text{traj err} = \norm{\uu-\uu^{\text{ref}}} = \sqrt{\norm{\pp-\pp^{\text{ref}}}^2 + \norm{\qq-\qq^{\text{ref}}}^2},
\end{align}
and the energy error 
\begin{align}
    \text{energy err} = \frac{\abs{H(\pp,\qq) - H(\pp^{\text{ref}}, \qq^{\text{ref}})}}{\abs{H\lr{\pp^{\text{ref}}, \qq^{\text{ref}}}}}. 
\end{align}

\subsection{Definition of the energy transform}

We first present our definition of the energy transform $\Lambda$, since both the Procrustes parareal and the NN solver will rely on this function. Based on the Hamiltonian \eqref{eq:H_fpu}, we define 
\begin{align}
    \Lambda_1: \; & \pp \in \mathbb{R}^{2m} \mapsto \frac{\pp}{\sqrt{2}} \in \mathbb{R}^{2m}, \\
    \Lambda_2: \; & \qq \in \mathbb{R}^{2m} \mapsto
    \begin{bmatrix}
        \mathbf{dq}_{\text{stiff}} \\
        \mathbf{dq}_{\text{soft}}
    \end{bmatrix} \in \mathbb{R}^{2m+1},
\end{align}
where
\begin{align}
    \mathbf{dq}_{\text{stiff}} :=
    \begin{bmatrix}
        \frac{\omega}{2} \lr{q_2 - q_1} \\[-0.5em] 
        \vdots \\[-0.5em]
        \frac{\omega}{2} \lr{q_{2i} - q_{2i-1}} \\[-0.5em]  
        \vdots \\[-0.5em] 
        \frac{\omega}{2} \lr{q_{2m} - q_{2m-1}}
    \end{bmatrix} \in \mathbb{R}^{m}, \quad 
    \mathbf{dq}_{\text{soft}} :=
    \begin{bmatrix}
        \lr{q_1 - q_0}^2 \\[-0.5em]  
        \vdots \\[-0.5em] 
        \lr{q_{2i+1} - q_{2i}}^2 \\[-0.5em]  
        \vdots \\[-0.5em] 
        \lr{q_{2m+1} - q_{2m}}^2 \\ 
    \end{bmatrix} \in \mathbb{R}^{m+1}.
\end{align}
Then $\Lambda$ can be written as 
\begin{align}
    \Lambda: \begin{bmatrix} \pp \\ \qq\end{bmatrix} \in \mathbb{R}^{4m} \mapsto 
    \begin{bmatrix}
        \Lambda_1(\pp) \\ 
        \Lambda_2(\qq) 
    \end{bmatrix} \in \mathbb{R}^{4m+1}.
\end{align}
One can check the $l_2$ norm squared of $\Lambda\lr{[\pp, \qq]}$ recovers \eqref{eq:H_fpu}.

To define the pseudo-inverse $\Lambda^{\dagger}$, the main task is to recover $\qq$ from a given vector $\Lambda_2(\qq)$. This can be done by solving a nonlinear least squares problem: given $\Tilde{\mathbf{dq}}_{\text{stiff}}, \Tilde{\mathbf{dq}}_{\text{soft}}$, find
\begin{align}
    \qq_* = \argmin_{\qq} \norm{\Lambda_2\lr{\qq} - 
    \begin{bmatrix}
        \Tilde{\mathbf{dq}}_{\text{stiff}} \\
        \Tilde{\mathbf{dq}}_{\text{soft}}
    \end{bmatrix}}_2^2.
\end{align}
We use an adaptive nonlinear least squares algorithm, called NL2SOL \cite{dennis1981algorithm}, to solve the problem. In the Procrustes parareal setup, $\Lambda^{\dagger}$ is only evaluated when applying the correction operator $\Psi^{(k)}_{\Delta t} := \Lambda^{\dagger} \Omega^{(k)}_{\Delta t}  \Lambda $ to some solution $\uu$. Since the correction by the unitary matrix $ \Omega^{(k)}_{\Delta t}$ is expected to be small, we use the $\qq$ component of $\uu$ as an initial guess in the least squares algorithm.

\subsection{The Procrustes Parareal Method}

In this section we discuss our choice of numerical integrators and present results for the Procrustes parareal method.

\subsubsection*{Choice of numerical solvers}

For Hamiltonian systems, we use symplectic integrators for the fine and coarse solvers. The stepsize $h$ of the integrator should be order $\mathcal{O}(\omega^{-1})$ or smaller for accurate integration. For the coarse solver $C_{\Delta t}$, we use a \nth{4}-order symplectic algorithm developed by Calvo and Sanz-Serna \cite{calvo1993development}. We consider two stepsizes $h=2^{-9}$ and $h=2^{-8}$ for comparison. For the fine solver $F_{\Delta t}$, we use an \nth{8}-order symplectic algorithm developed by Kahan and Li \cite{kahan1997composition} with a stepsize $h=2^{-18}$. The numerical integrators are denoted by $\Phi^{\text{CSS4}, h=2^{-9}}_{\Delta t}$, $\Phi^{\text{CSS4}, h=2^{-8}}_{\Delta t}$ and $\Phi^{\text{KL8}, h=2^{-18}}_{\Delta t}$.

Because we aim to simulate for a long time interval~$\mathcal{O}(\omega)$, generating the reference trajectory requires $\omega/h\approx 8\times10^7$ fine steps. Hence, we perform fine solver computations in quadruple precision to avoid loss of significant digits. The computations of the coarse solver are done in double precision. We perform all numerical computations in Julia. We use the MultiFloats library to obtain quadruple precision numbers.  

\begin{figure}[h!]
    \centering
    \begin{subfigure}[b]{0.7\linewidth}
        \centering
        \includegraphics[width=\linewidth]{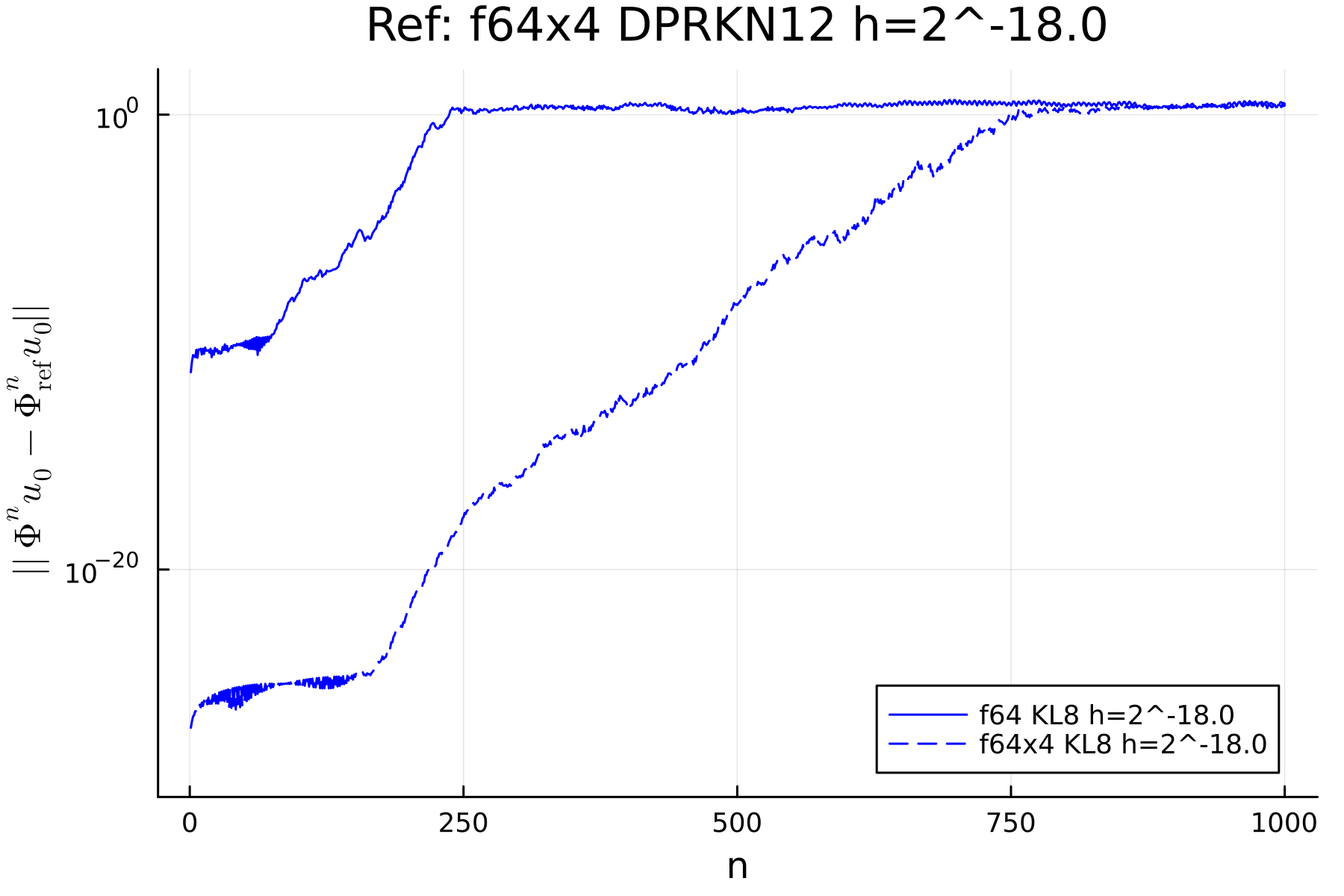}
        \caption{trajectory error}
    \end{subfigure}
    \hfill
    \begin{subfigure}[b]{0.7\linewidth}
        \centering
        \includegraphics[width=\linewidth]{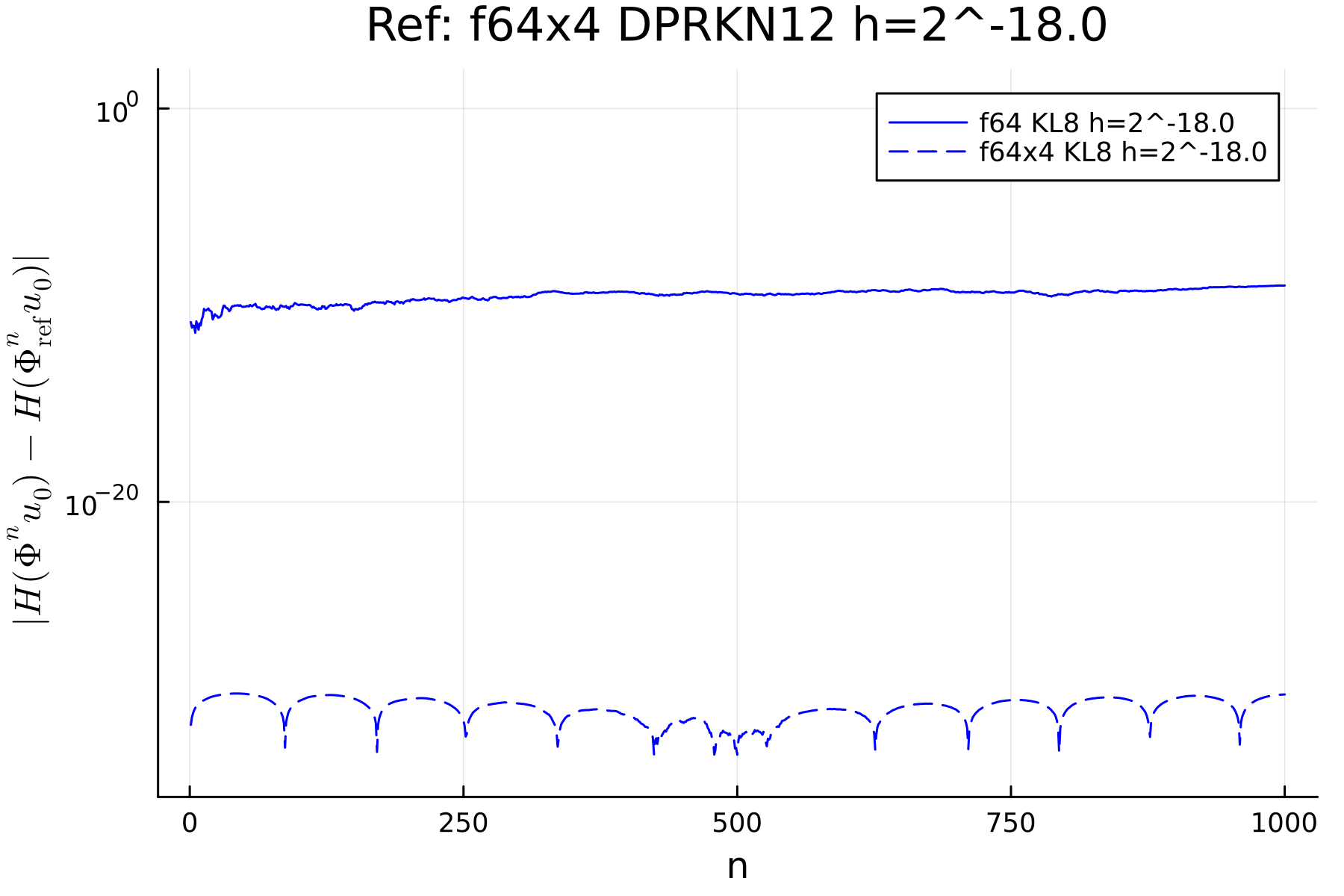}
        \caption{energy error}
    \end{subfigure}
    \caption{Errors in trajectories generated by applying the fine solver $F_{\Delta t}=\Phi^{\text{KL8}, h=2^{-18}}_{\Delta t}$, implemented in double precision and quadruple precision, sequentially for $N=1000$ steps of $\Delta t=1.0$. F64 refers to double precision and f64x4 refers to quadruple precision. The numerical integrator used for generating the reference trajectory here is $\Phi^{\text{DPRKN12}, h=2^{-18}}_{\Delta t}$, implemented in quadruple precision.}
    \label{fig:f64_vs_f64x4_kl8}
\end{figure}

Figure~\ref{fig:f64_vs_f64x4_kl8} shows the global errors of the fine solutions computed up to $T=1000$ in double precision versus in quadruple precision. Here we use a method of even higher order, the \nth{12}-order explicit Runge-Kutta-Nyström method, with stepsize $h^{-18}$ implemented in quadruple precision to serve as the reference solution map. We can see the trajectory errors grow over time while the energy errors are stable (expected since it is a symplectic integrator) for both precision. The trajectory errors grow linearly in time at first, and then grow exponentially. Based on the trajectory errors, we conclude the fine solutions computed in quadruple precision lose digits much later than fine solutions computed in double precision. Even with quadruple precision, the fine solutions are not reliable after $n=500$. In the rest of the paper, unless otherwise mentioned, we compare the trajectory errors against the reference only up to $n=500$. For energy errors, we may compare for $n>500$ since the reference energy is almost a constant. 

\begin{figure}[h!]
    \centering
    \begin{tabular}{cccc}
        &  Traj error & Energy error \\
        \small{\makecell{$C_{\Delta t}$ \\ quadruple \\ precision }} & \adjincludegraphics[height=4cm, valign=m, trim={{0.04\width} 0 {0.06\width} 0}, clip]{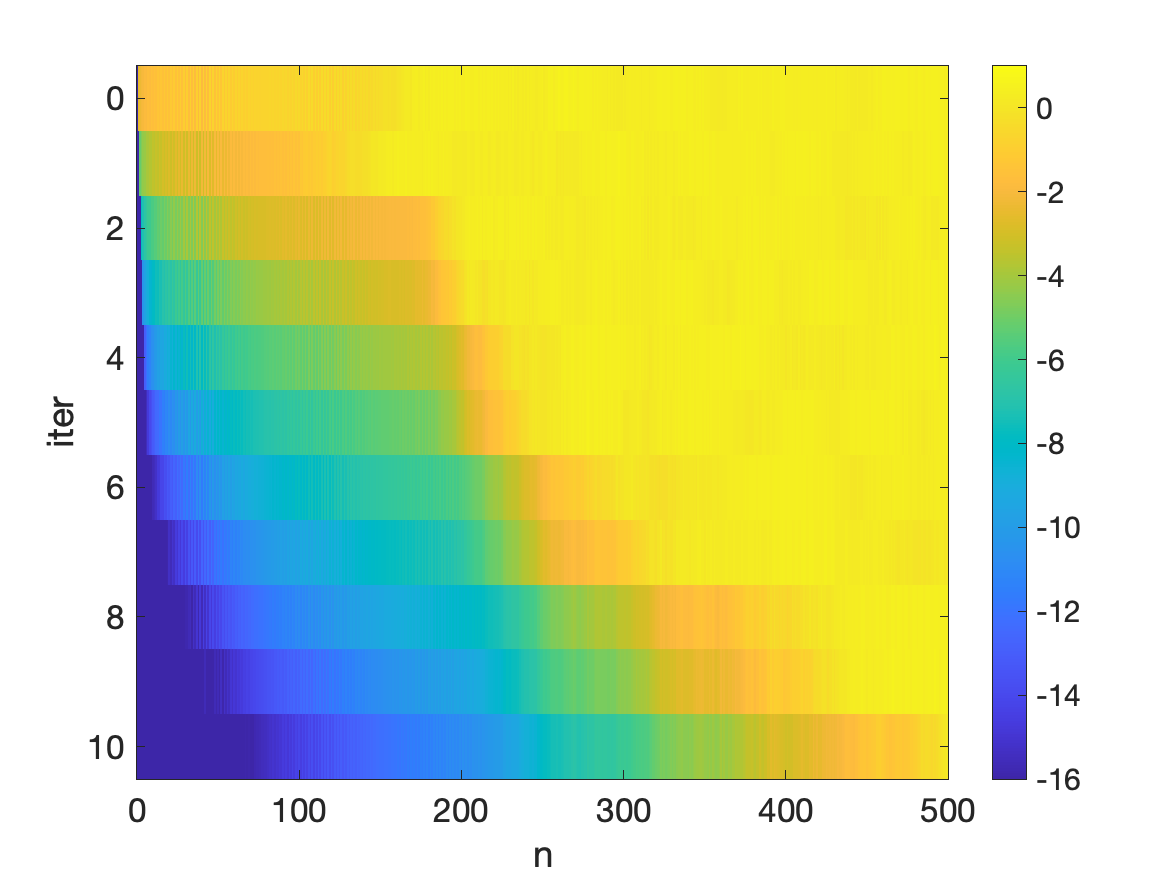} & \adjincludegraphics[height=4cm, valign=m, trim={{0.04\width} 0 {0.06\width} 0}, clip]{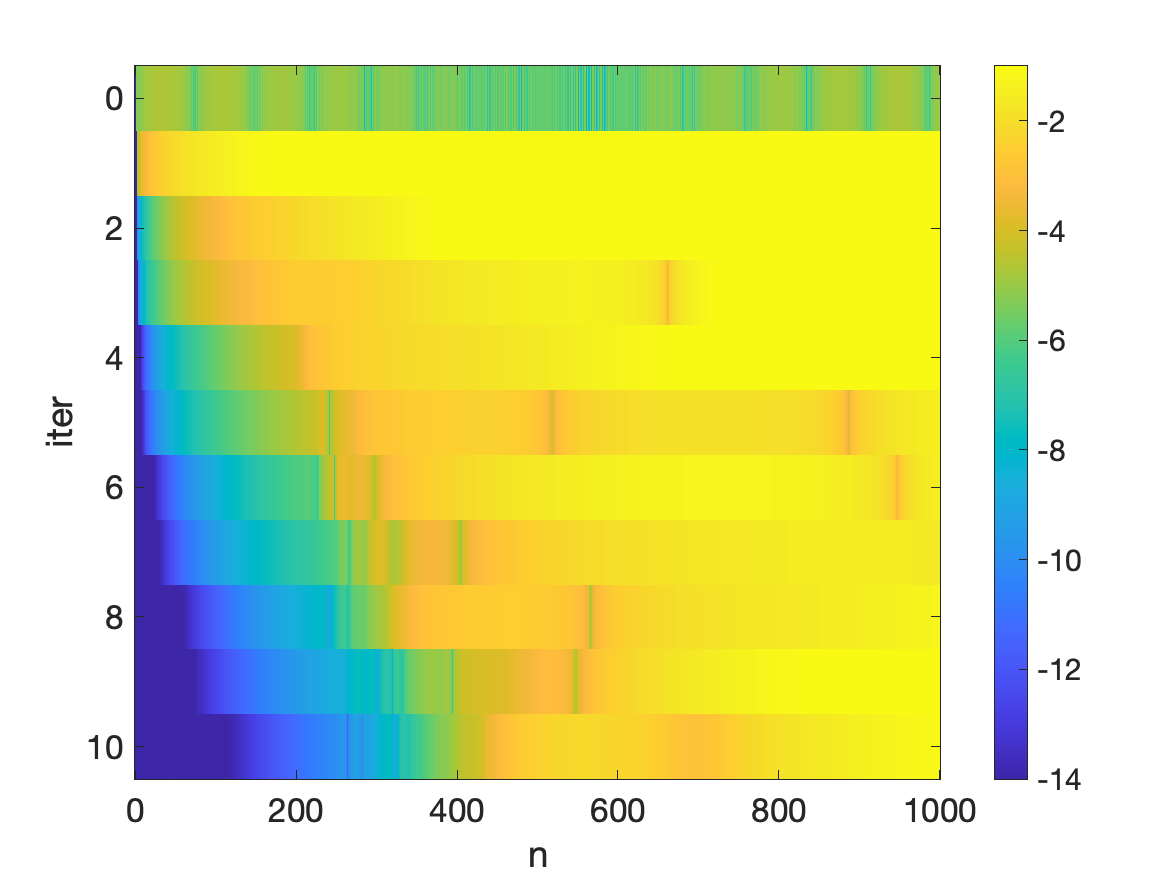} \\
        \small{\makecell{$C_{\Delta t}$ \\ double \\ precision  }} &  \adjincludegraphics[height=4cm, valign=m, trim={{0.04\width} 0 {0.06\width} 0}, clip]{figs/parareal/plain/css4-h=2e-9/log10_traj_err_N500_K10} & \adjincludegraphics[height=4cm, valign=m, trim={{0.04\width} 0 {0.06\width} 0}, clip]{figs/parareal/plain/css4-h=2e-9/log10_H_err_N1000_K10} 
    \end{tabular}
    \caption{Log (base 10) errors in plain parareal solutions computed with $C_{\Delta t} =\Phi^{\text{CSS4}, h=2^{-9}}_{\Delta t}$ implemented in double precision and quadruple precision, and $F_{\Delta t}=\Phi^{\text{KL8}, h=2^{-18}}_{\Delta t}$ implemented in quadruple precision only ($\Delta t=1.0$, $T=1000$).}
    \label{fig:f64_vs_f64x4_parareal}
\end{figure}

We also observe an effect of floating point precision on parareal iterations. Figure~\ref{fig:f64_vs_f64x4_parareal} shows the errors in parareal solutions computed with a coarse solver implemented in double precision versus quadruple precision. The fine solver is fixed using quadruple precision. We found the error plots differ significantly after $n=200$ steps. Using double precision for the coarse solver prevents the parareal solutions from improving after $n=200$. Unfortunately, we have to use double precision for the coarse solver for the rest of comparisons, given the facts that 1. the library for the least squares algorithm involved in inverting $\Lambda$ only supports double precision, and 2. the NN solver is double precision.

\subsubsection*{Plain versus Procrustes}

Using the introduced numerical solvers, we run the plain parareal method and the Procrustes parareal method for $N=1000$ steps and $k=10$ iterations, from the same initial condition. 

\begin{figure}[h!]
    \centering
    \begin{tabular}{cccc}
        &  $k=0$ & $k=3$ \\
        Plain & 
        \adjincludegraphics[height=4cm, valign=m, trim={{0.04\width} 0 {0.06\width} 0}, clip]{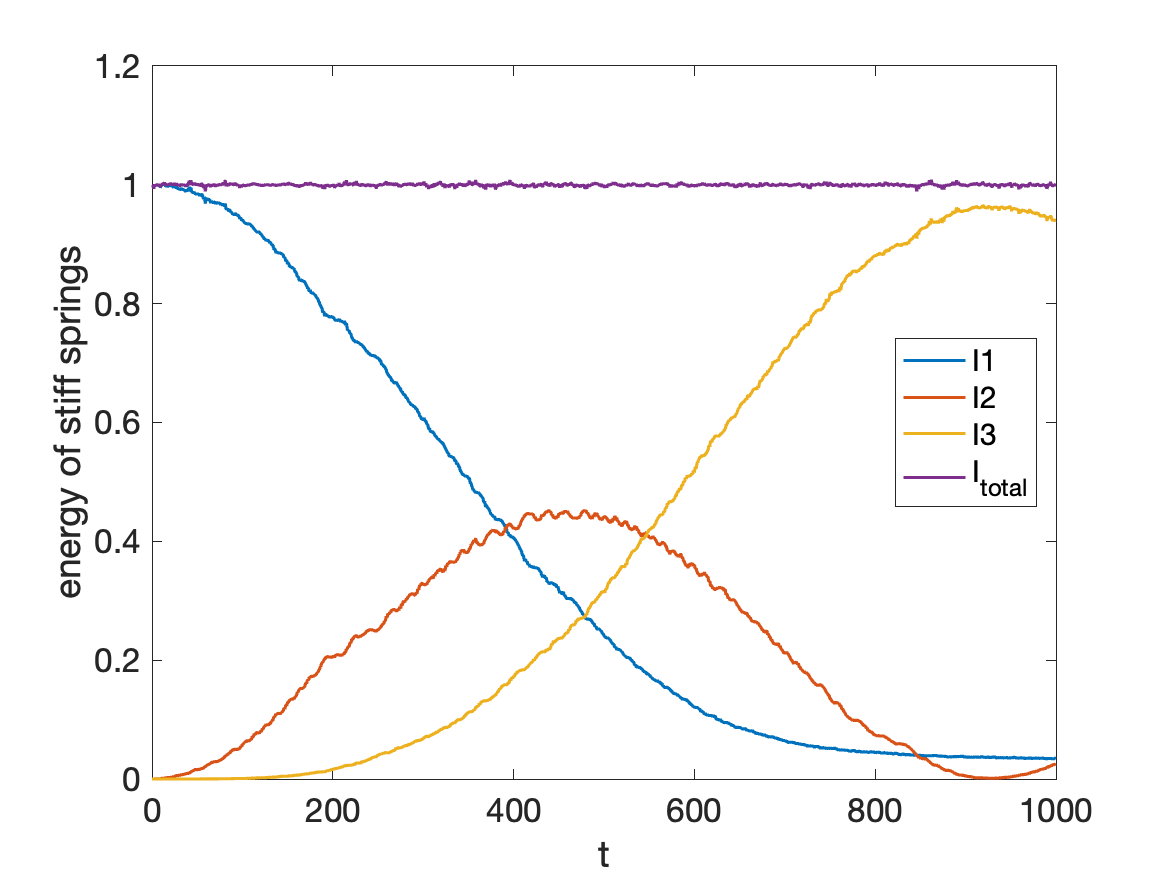} & \adjincludegraphics[height=4cm, valign=m, trim={{0.04\width} 0 {0.06\width} 0}, clip]{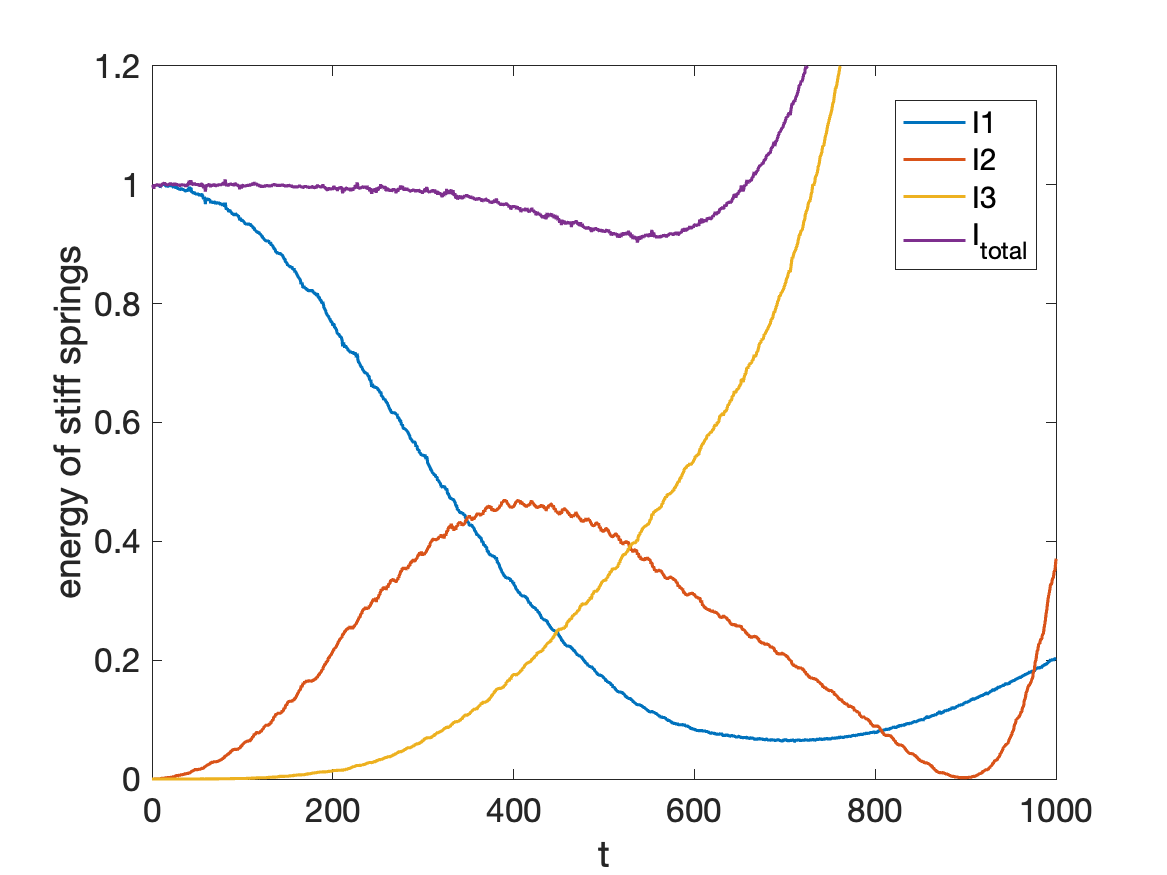} 
         \\
        Procrustes & \adjincludegraphics[height=4cm, valign=m, trim={{0.04\width} 0 {0.06\width} 0}, clip]{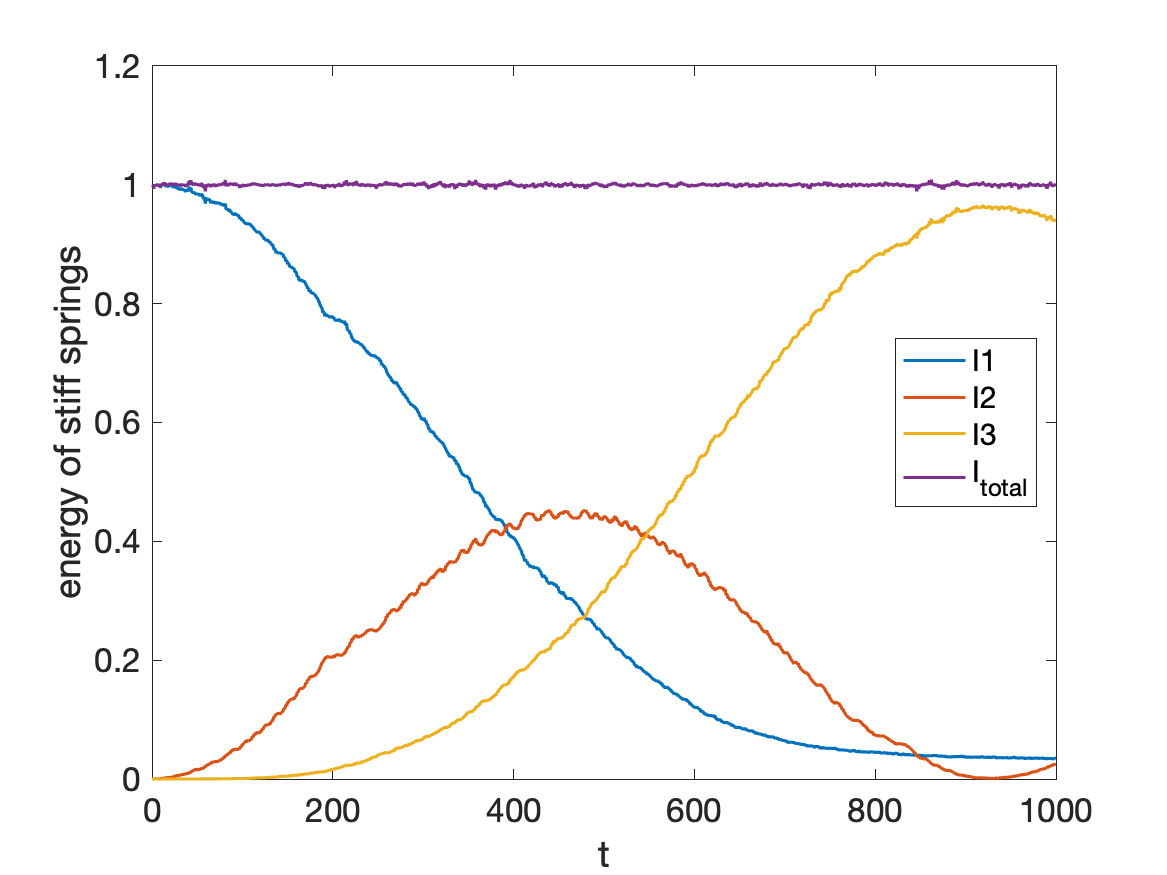} & \adjincludegraphics[height=4cm, valign=m, trim={{0.04\width} 0 {0.06\width} 0}, clip]{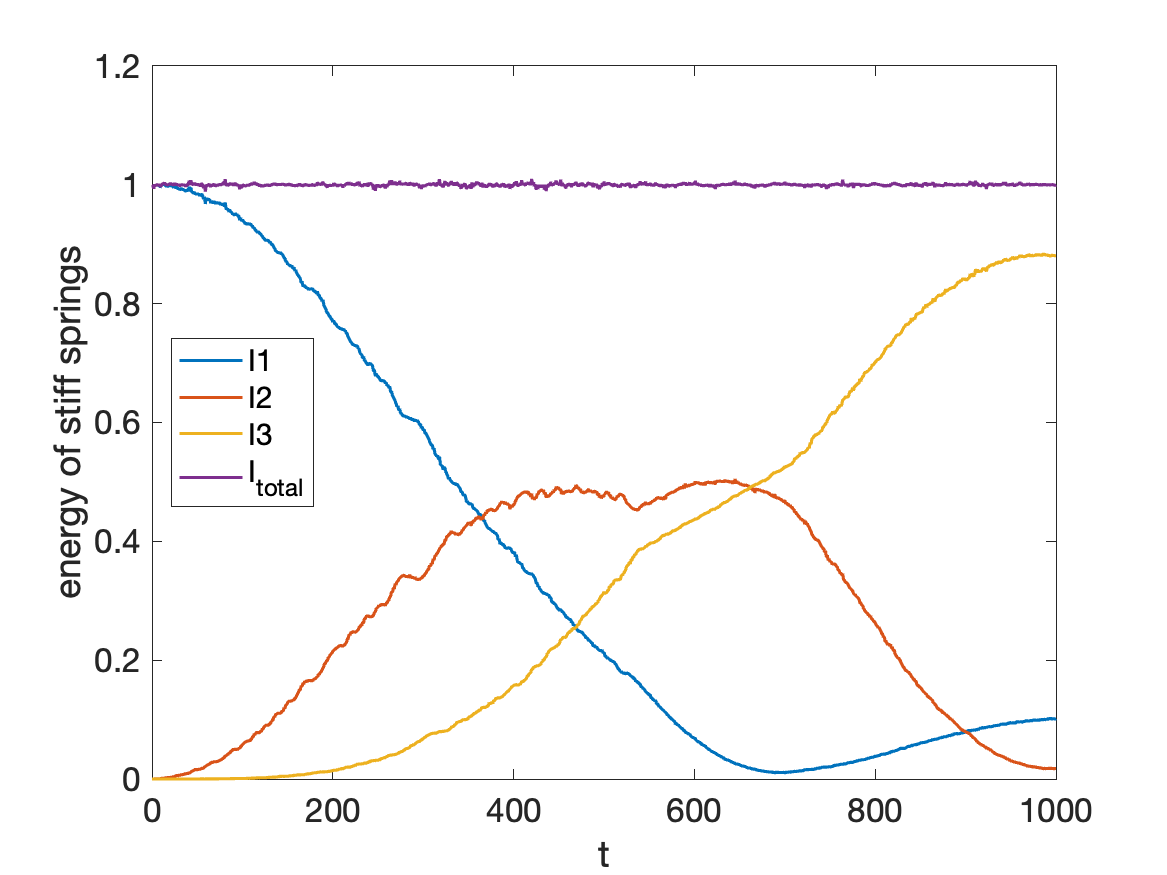} 
    \end{tabular}
    \caption{Energy profiles of the stiff springs computed from plain parareal solutions and Procrustes parareal solutions at iteration~0 and iteration~3 ($\Delta t=1.0$, $T=1000$,  $C_{\Delta t} =\Phi^{\text{CSS4}, h=2^{-9}}_{\Delta t}$,  $F_{\Delta t}=\Phi^{\text{KL8}, h=2^{-18}}_{\Delta t}$ ).}
    \label{fig:parareal-stiffenergies-css4-h=2e-9}
\end{figure}

\begin{figure}[h!]
    \centering
    \includegraphics[width=0.7\linewidth]{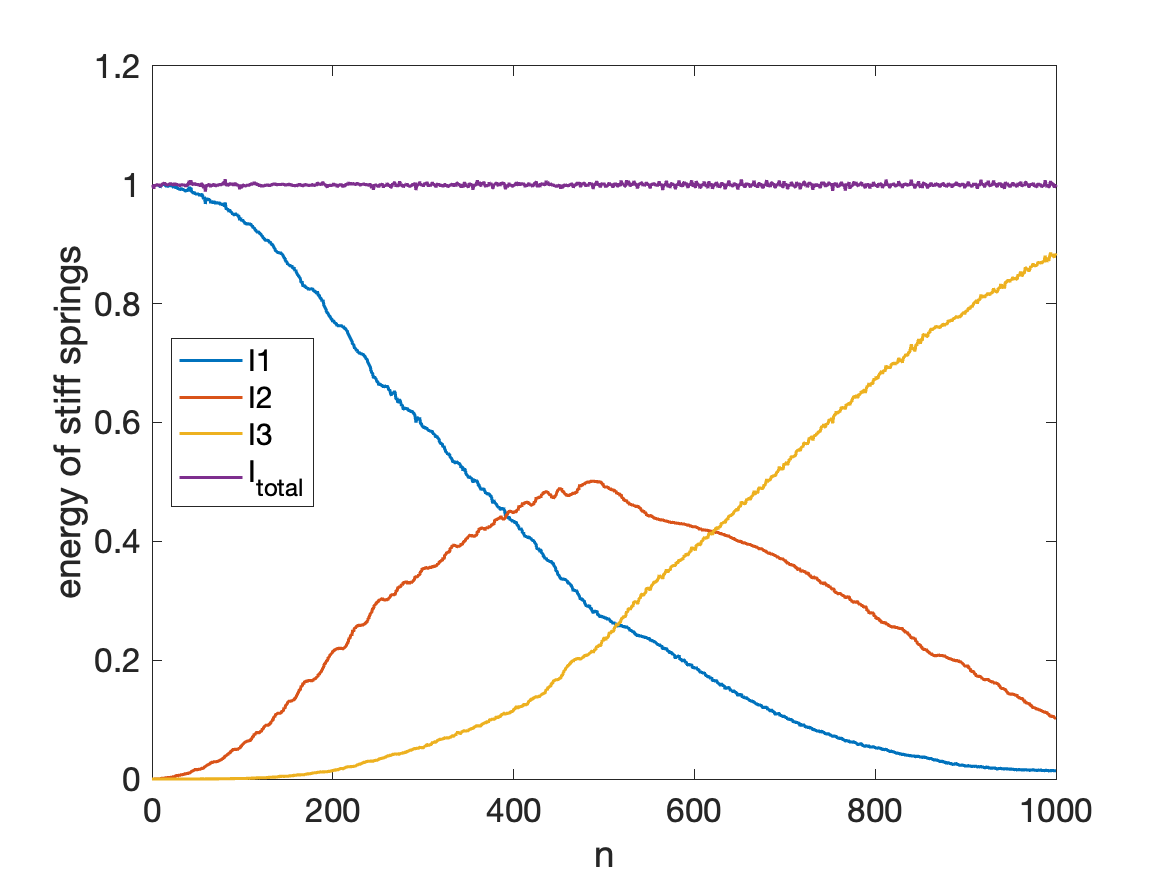}
    \caption{Reference energy profile of the stiff springs ($\Delta t=1.0$, $T=1000$).}
    \label{fig:ref-stiffenergies}
\end{figure}

Figure~\ref{fig:parareal-css4-h=2e-9} shows errors in the computed trajectories for $C_{\Delta t} =\Phi^{\text{CSS4}, h=2^{-9}}_{\Delta t}$. Overall, both methods improve the initial coarse solution, but the improvement becomes small after several iterations. Comparing the trajectory error plots, we observe the Procrustes parareal solutions improve faster than the plain parareal solutions over iterations. Comparing the energy error plots, we see the Procrustes parareal solutions not only improve faster, they are also more stable in energy than the plain parareal solutions (in particular, see how the energy errors grow from iteration~0 to iteration~1 in plain parareal versus in Procrustes parareal). The stability issue can be seen more clearly in Figure~\ref{fig:parareal-stiffenergies-css4-h=2e-9}, where we plot energy of the three stiff springs and their total energy from the computed solutions. As shown in the reference energy profile in Figure~\ref{fig:ref-stiffenergies}, during the simulated time range, there is energy exchange among the stiff springs while the total energy of the stiff springs remains almost a constant. For the plain parareal solution at iteration~3, the total energy blows up by the end of the simulation. On the contrary, the total energy for the Procrustes parareal solution is almost conserved. 

\begin{figure}[h!]
    \centering
    \begin{tabular}{cccc}
        &  $k=0$ & $k=3$ \\
        Plain & 
        \adjincludegraphics[height=4cm, valign=m, trim={{0.04\width} 0 {0.06\width} 0}, clip]{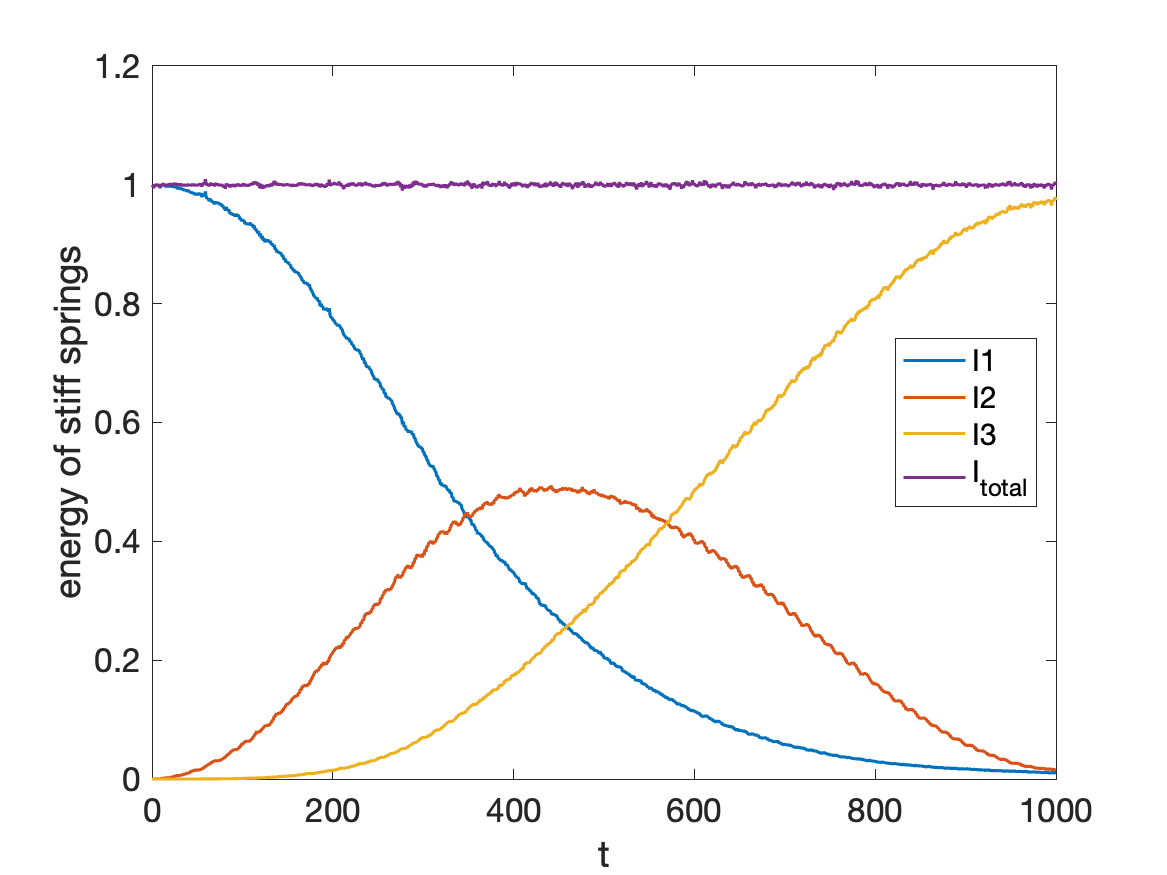} & \adjincludegraphics[height=4cm, valign=m, trim={{0.04\width} 0 {0.06\width} 0}, clip]{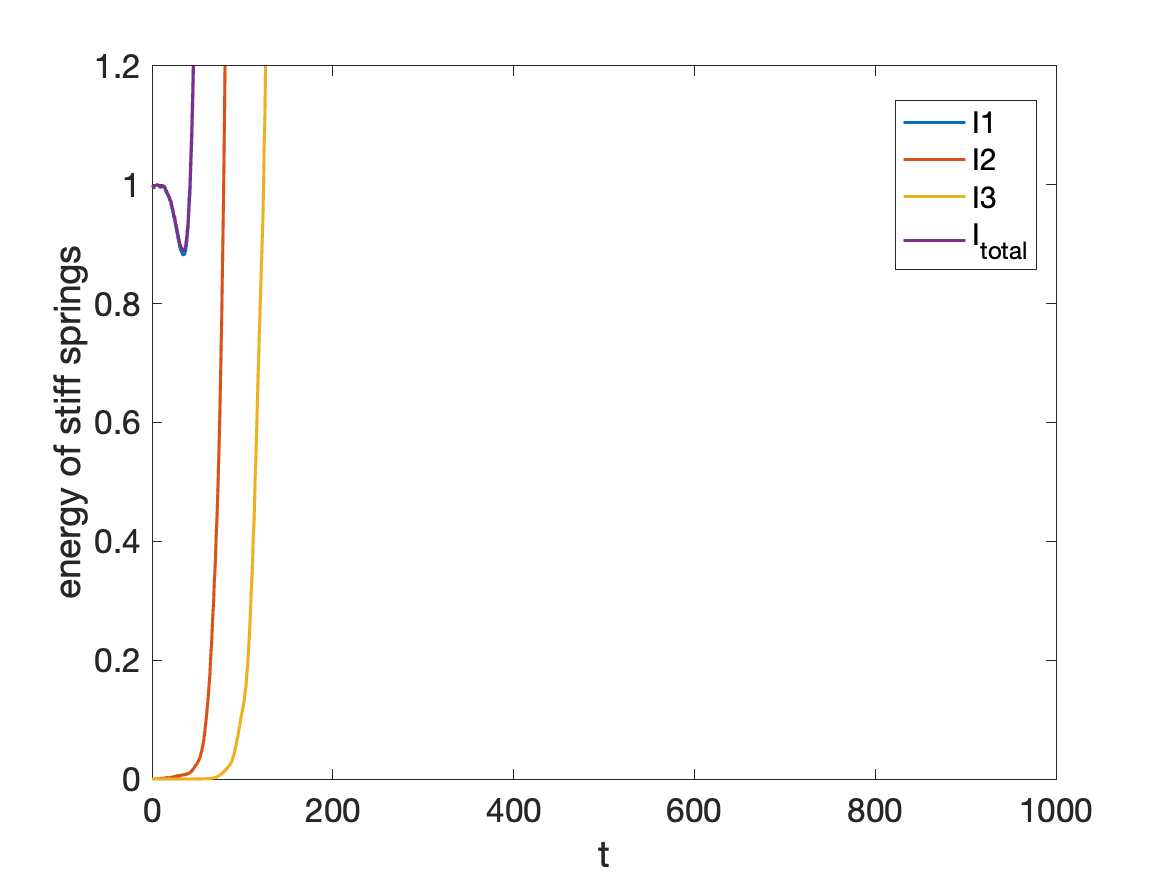} 
         \\
        Procrustes & \adjincludegraphics[height=4cm, valign=m, trim={{0.04\width} 0 {0.06\width} 0}, clip]{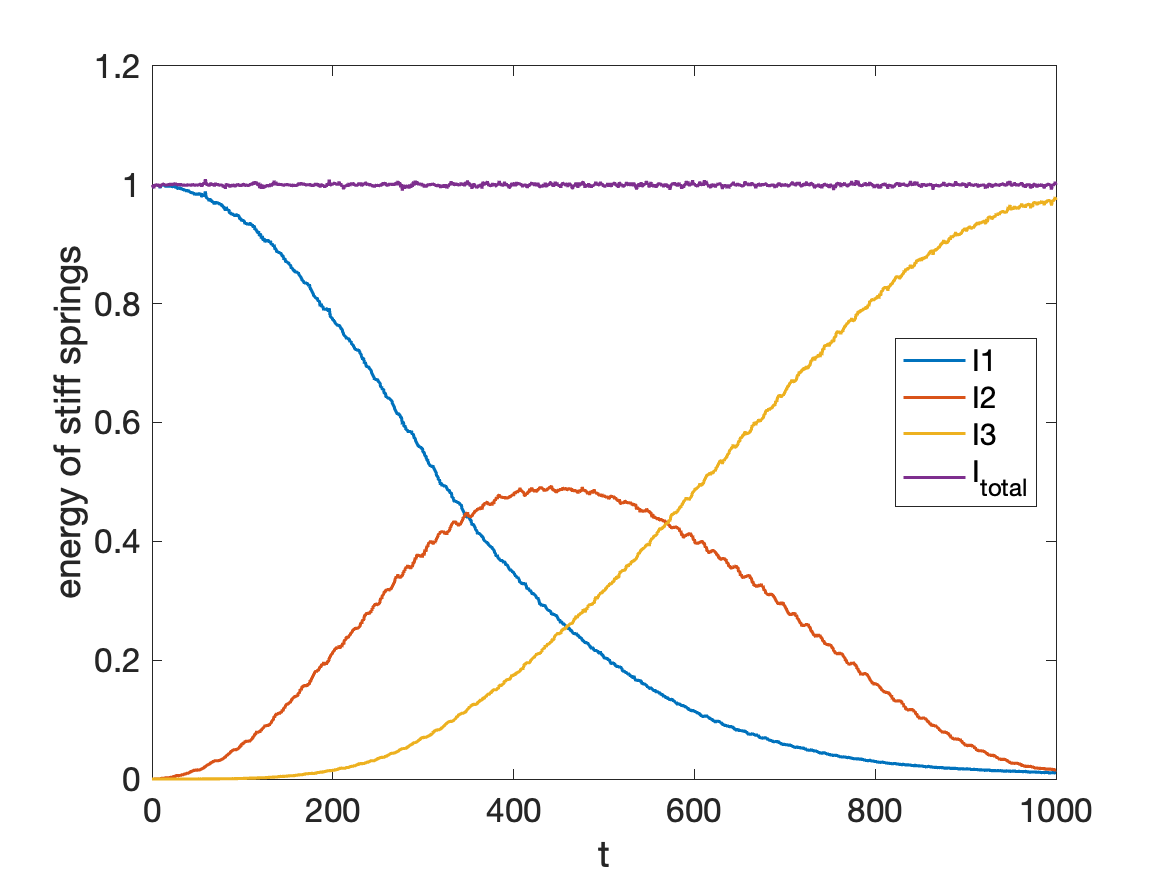} & \adjincludegraphics[height=4cm, valign=m, trim={{0.04\width} 0 {0.06\width} 0}, clip]{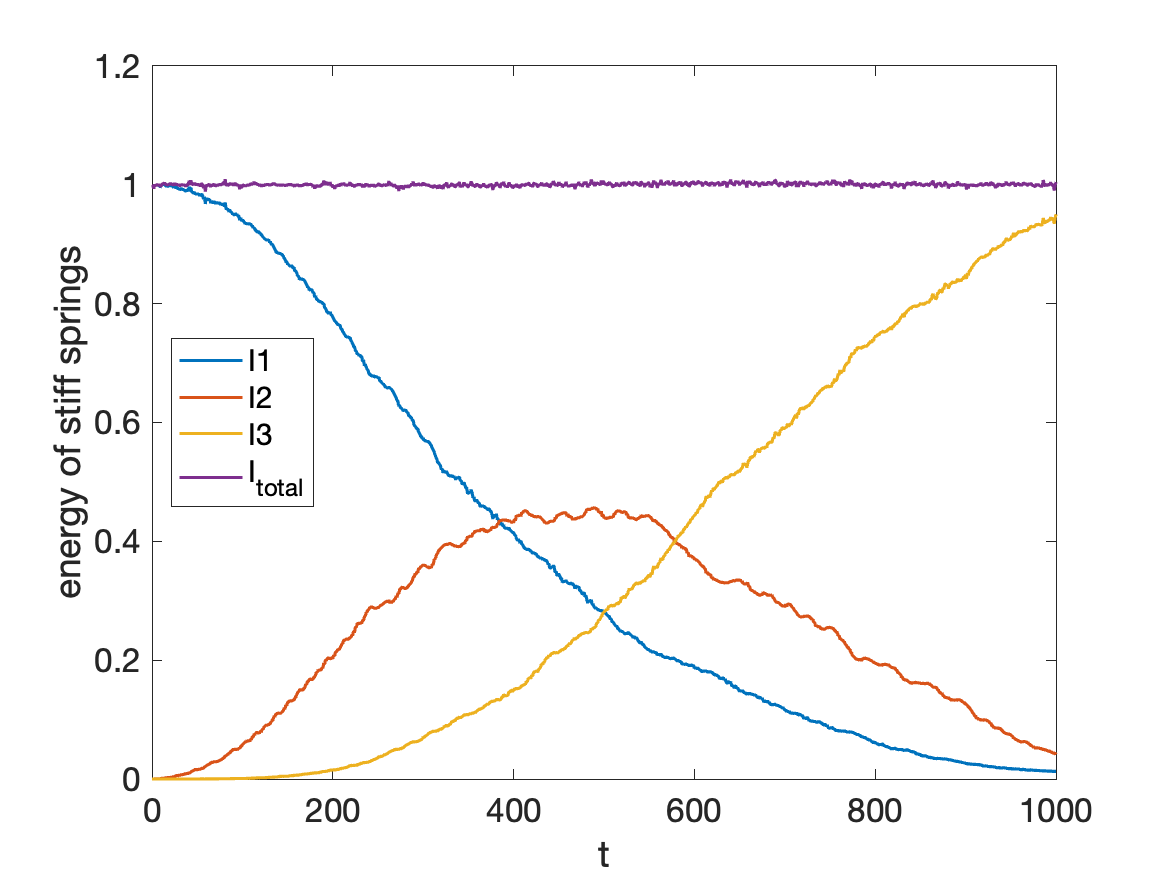} 
    \end{tabular}
    \caption{Energy profiles of the stiff springs computed from plain parareal solutions and Procrustes parareal solutions at iteration~0 and iteration~3 ($\Delta t=1.0$, $T=1000$,  $C_{\Delta t} =\Phi^{\text{CSS4}, h=2^{-8}}_{\Delta t}$,  $F_{\Delta t}=\Phi^{\text{KL8}, h=2^{-18}}_{\Delta t}$ ).}
    \label{fig:parareal-stiffenergies-css4-h=2e-8}
\end{figure}

We repeat the comparison for a less accurate coarse solver $C_{\Delta t} =\Phi^{\text{CSS4}, h=2^{-8}}_{\Delta t}$ in Figure~\ref{fig:parareal-css4-h=2e-8} and Figure~\ref{fig:parareal-stiffenergies-css4-h=2e-8}. Comparing Figure~\ref{fig:parareal-css4-h=2e-8} and Figure~\ref{fig:parareal-css4-h=2e-9}, we see the improvement of Procrustes parareal over plain parareal is more substantial for the less accurate coarse solver.

\begin{figure}[h!]
    \centering
    \begin{tabular}{cccc}
        &  Traj error &  Energy error \\
        Plain  & \adjincludegraphics[height=4cm, valign=m, trim={{0.04\width} 0 {0.06\width} 0}, clip]{figs/parareal/plain/css4-h=2e-9/log10_traj_err_N500_K10} & \adjincludegraphics[height=4cm, valign=m, trim={{0.04\width} 0 {0.06\width} 0}, clip]{figs/parareal/plain/css4-h=2e-9/log10_H_err_N1000_K10} \\
        Procrustes  &  \adjincludegraphics[height=4cm, valign=m, trim={{0.04\width} 0 {0.06\width} 0}, clip]{figs/parareal/phasecorr/css4-h=2e-9/log10_traj_err_N500_K10} & \adjincludegraphics[height=4cm, valign=m, trim={{0.04\width} 0 {0.06\width} 0}, clip]{figs/parareal/phasecorr/css4-h=2e-9/log10_H_err_N1000_K10} 
    \end{tabular}
    \caption{Log (base 10) errors in plain parareal solutions and Procrustes parareal solutions ($\Delta t=1.0$, $T=1000$, $C_{\Delta t} =\Phi^{\text{CSS4}, h=2^{-9}}_{\Delta t}$,  $F_{\Delta t}=\Phi^{\text{KL8}, h=2^{-18}}_{\Delta t}$ ).}
    \label{fig:parareal-css4-h=2e-9}
\end{figure}

\begin{figure}[h!]
    \centering
    \begin{tabular}{cccc}
        &  Traj error &  Energy error \\
        Plain  & \adjincludegraphics[height=4cm, valign=m, trim={{0.04\width} 0 {0.06\width} 0}, clip]{figs/parareal/plain/css4-h=2e-8/log10_traj_err_N500_K10} & \adjincludegraphics[height=4cm, valign=m, trim={{0.04\width} 0 {0.06\width} 0}, clip]{figs/parareal/plain/css4-h=2e-8/log10_H_err_N1000_K10} \\
        Procrustes  &  \adjincludegraphics[height=4cm, valign=m, trim={{0.04\width} 0 {0.06\width} 0}, clip]{figs/parareal/phasecorr/css4-h=2e-8/log10_traj_err_N500_K10} & \adjincludegraphics[height=4cm, valign=m, trim={{0.04\width} 0 {0.06\width} 0}, clip]{figs/parareal/phasecorr/css4-h=2e-8/log10_H_err_N1000_K10} 
    \end{tabular}
    \caption{Log (base 10) errors in plain parareal solutions and Procrustes parareal solutions ($\Delta t=1.0$, $T=1000$, $C_{\Delta t} =\Phi^{\text{CSS4}, h=2^{-8}}_{\Delta t}$,  $F_{\Delta t}=\Phi^{\text{KL8}, h=2^{-18}}_{\Delta t}$ ).}
    \label{fig:parareal-css4-h=2e-8}
\end{figure}

\subsection{NN solution map}

In this section, we present the setups for learning the solution map $\Phi^{\text{NN}}_{\Delta t}$ and study how its quality is affected by different options of training data, loss function, and network architecture. To evaluate performance of a learned $\Phi^{\text{NN}}_{\Delta t}$, we generate a 1000-step trajectory from the initial condition \eqref{eq:test_u0} by sequential applications of $\Phi^{\text{NN}}_{\Delta t}$.

Using the proposed data generation algorithms, we generated two sets of input data, denoted by $\mathcal{D}_0^{\text{HMC-$H_0$}}$ and $\mathcal{D}_0^{\text{TrajEnsemble-$H_0$}}$ where $H_0 = H\lr{\pp_{\text{init}},\qq_{\text{init}}}$. The parameters for data generation algorithms are given in Table~\ref{tab:training_data_input}. Each dataset has 200k examples of inputs $\uu_0$. We emphasize that $\lr{\pp_{\text{init}},\qq_{\text{init}}}$ is not included in either dataset even though $\qq_{\text{init}}$ was used in the parameters. This is because both algorithms randomly sample $\pp$ given $\qq_{\text{init}}$. 

Given a $\mathcal{D}_0$, we generate the full training dataset $\mathcal{D}$ by propagating each $\uu_0$ for 5 steps using the fine solver $F_{\Delta t}$, and then collecting the input and target sequence pairs:
\begin{align}
    \mathcal{D} = \left\{ \lr{\uu_0, \{ \lr{F_{\Delta t}}^{i} \uu_0\}_{i=1,\cdots,5}}: \uu_0 \in \mathcal{D}_0  \right\}.
\end{align}
This leads to two training datasets $\mathcal{D}^{\text{TrajEnsemble-$H_0$}}$ and  $\mathcal{D}^{\text{HMC-$H_0$}}$.

Let ResNet($L$, $n$) denote a network with $L$ hidden layers and $n$ nodes per hidden layer. We considered several ResNet architectures, including a shallow network ResNet(4, 1000) and a deep network ResNet(75, 200). The two networks have similar number of trainable parameters, which is around $3\times 10^6$. Performances of the two networks are similar. Hence we will just report results for ResNet(4, 1000). 

The neural networks are implemented using PyTorch. We trained the networks using the mini-batch Adam algorithm \cite{kingma2014adam} with weight decay \cite{loshchilov2017decoupled}. To accelerate the training procedure, we used a 1-cycle learning rate scheduler \cite{smith2019super} that anneals the learning rate from an initial value to some maximum value and then to some minimum value within a fixed number of epochs. We used 10000 epochs to train a ResNet(4, 1000) and 5000 epochs to train a ResNet(75, 200). 

\begin{table}[h!]
    \caption{Descriptions of input data sets.}
    \label{tab:training_data_input}
    \centering
    \begin{tabular}{lM}
        \toprule
        Dataset & Parameters \\ 
        \midrule
        $\mathcal{D}_0^{\text{TrajEnsemble-$H_0$}}$ & \makecell{$\qq_0 = \qq_{\text{init}}$, $\sigma=0.1$, $N_{\mathrm{level sets}}=400$, $N_{\mathrm{traj}}=10$, \\ $L=50$, $\delta t=0.1$, $F_{\delta t}=\Phi^{\text{CSS4}, h=5^{-6}}_{\delta t}$}  \\
        $\mathcal{D}_0^{\text{HMC-$H_0$}}$ & \makecell{$\qq_0 = \qq_{\text{init}}$, $\sigma=0.1$, $N_{\mathrm{chains}}=100$, $N_{\mathrm{trans}}=2000$, \\ $\delta t=0.4$, $F_{\delta t}=\Phi^{\text{CSS4}, h=5^{-6}}_{\delta t}$} \\
        \botrule 
    \end{tabular}
\end{table}

\subsubsection*{Effects of training data}

\begin{figure}[h!]
    \centering
    \begin{subfigure}[b]{0.48\linewidth}
        \centering
        \includegraphics[width=\linewidth]{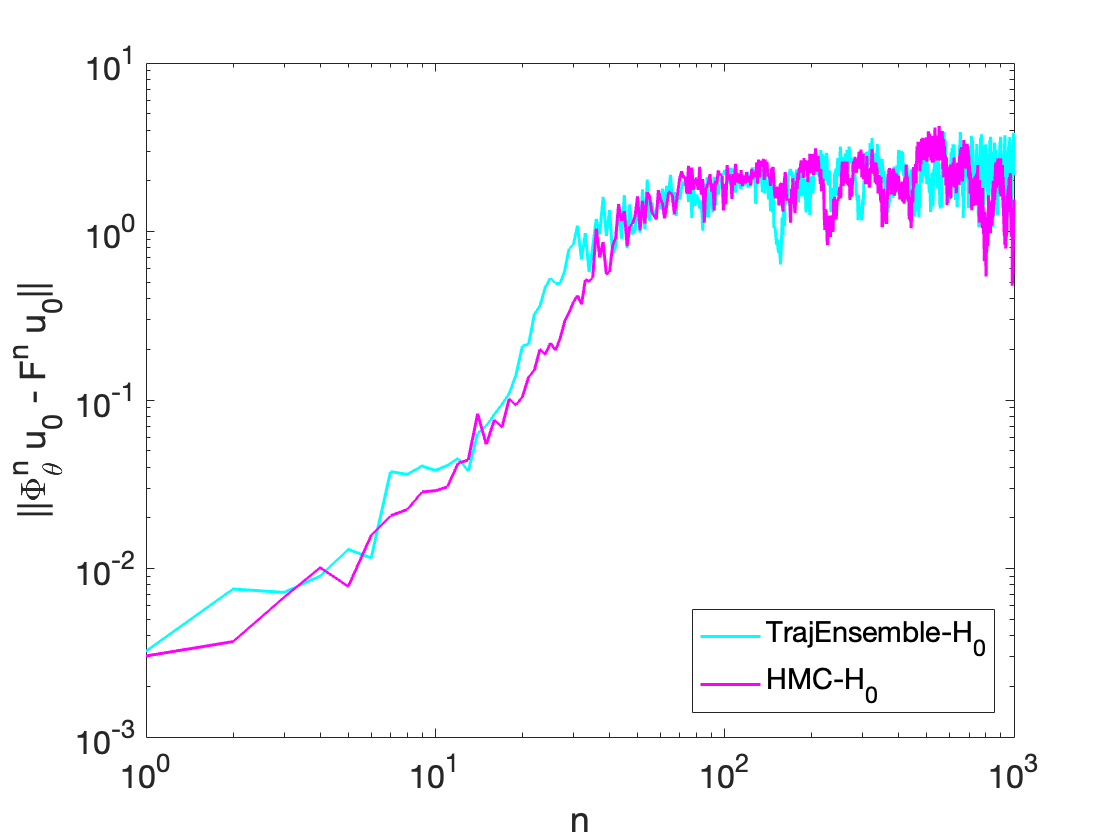}
        \caption{trajectory error}
    \end{subfigure}
    \hfill
    \begin{subfigure}[b]{0.48\linewidth}
        \centering
        \includegraphics[width=\linewidth]{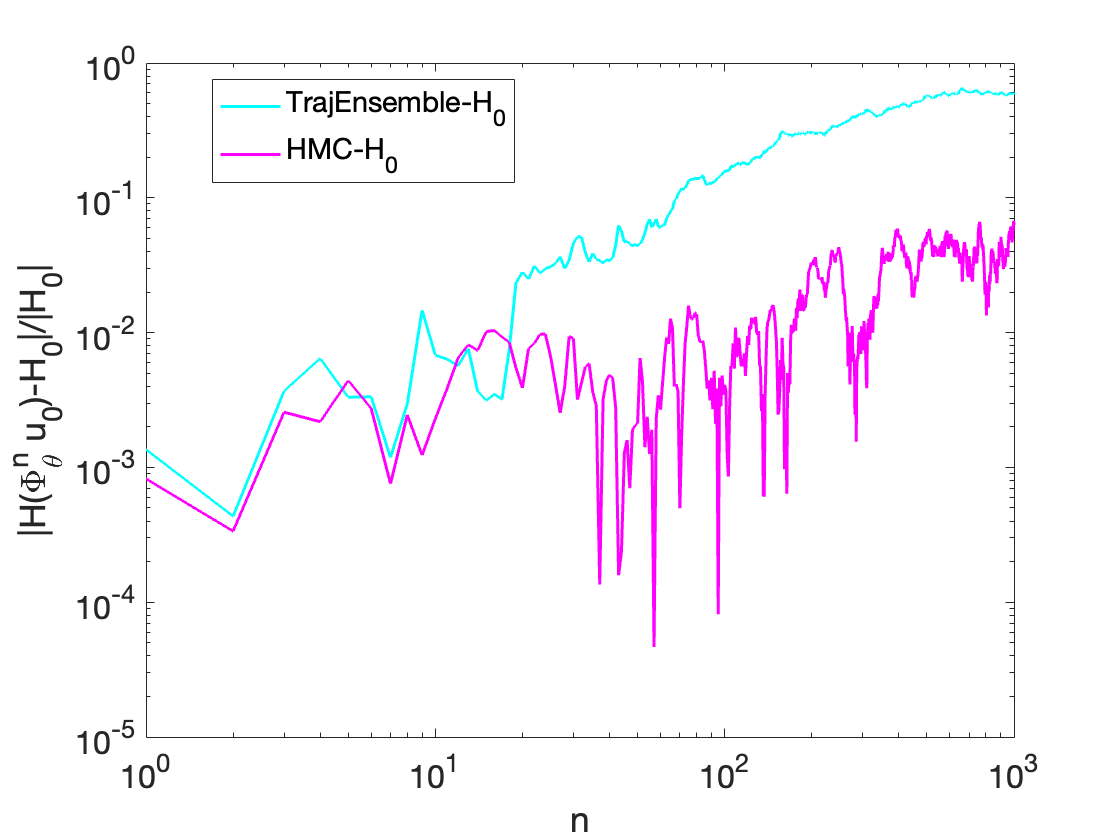}
        \caption{energy error}
    \end{subfigure}
    \caption{Errors in trajectories generated by NN solvers learned with different data.}
    \label{fig:effects_of_data}
\end{figure}

\begin{figure}[h!]
    \centering
    \includegraphics[width=0.7\linewidth]{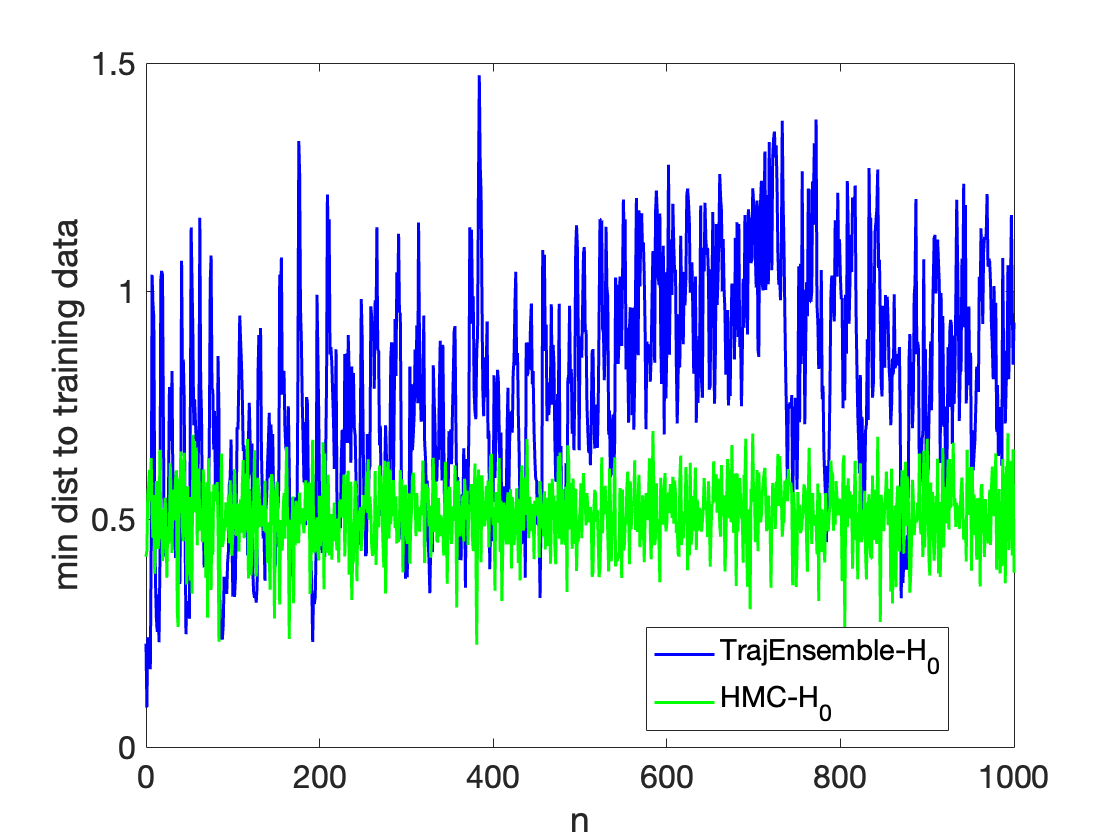}
    \caption{Minimum distance from each point on the reference trajectory to different training data in the phase space.}
    \label{fig:min_dist_to_data}
\end{figure}

We trained the shallow network ResNet(4, 1000) with different datasets, and with the one-step ($S=1$) MSE loss. Figure~\ref{fig:effects_of_data} shows the network trained with $\mathcal{D}^{\text{HMC-$H_0$}}$ is better than the network trained with $\mathcal{D}^{\text{TrajEnsemble-$H_0$}}$, especially in terms of the energy stability. We attribute this to the fact that the set of inputs $\mathcal{D}_0^{\text{HMC-$H_0$}}$ better represent the target distribution. As shown in Figure~\ref{fig:min_dist_to_data}, the minimum distance from each point along the reference trajectory to the set $\mathcal{D}_0^{\text{HMC-$H_0$}}$ is on average a lot smaller than the minimum distance to the set $\mathcal{D}_0^{\text{TrajEnsemble-$H_0$}}$. What is more, the minimum distance to $\mathcal{D}_0^{\text{TrajEnsemble-$H_0$}}$ increases along the trajectory, while the minimum distance to $\mathcal{D}_0^{\text{HMC-$H_0$}}$ stays stable over time. 

\subsubsection*{Effects of sequence length in loss function}

\begin{figure}[h!]
    \centering
    \begin{subfigure}[b]{0.48\linewidth}
        \centering
        \includegraphics[width=\linewidth]{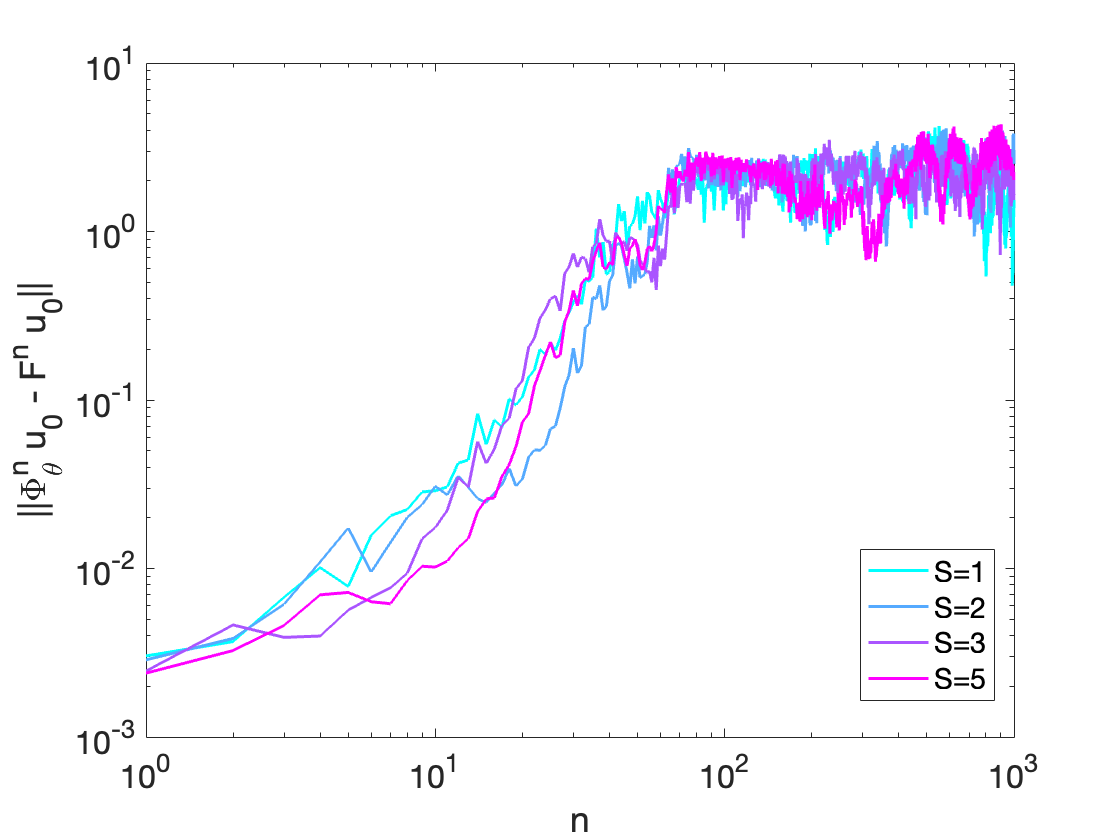}
        \caption{trajectory error}
    \end{subfigure}
    \hfill
    \begin{subfigure}[b]{0.48\linewidth}
        \centering
        \includegraphics[width=\linewidth]{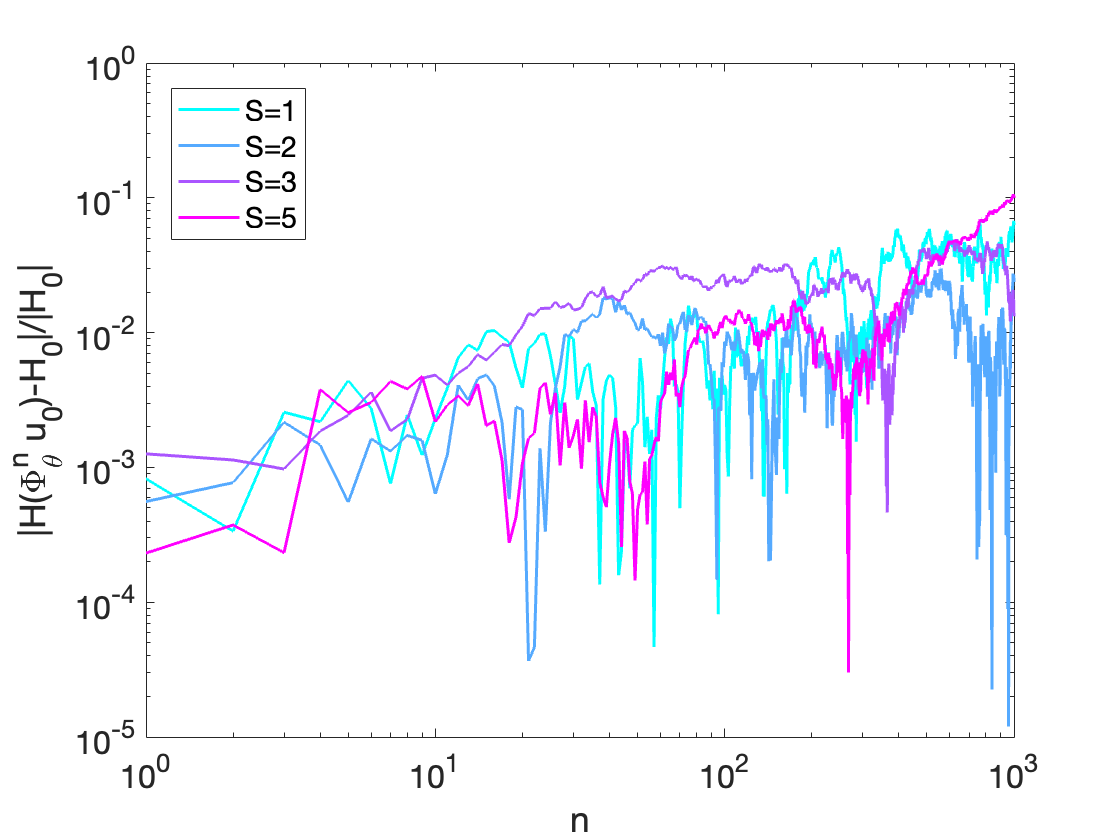}
        \caption{energy error}
    \end{subfigure}
    \caption{Errors in trajectories generated by NN solvers learned with different sequence length $S$ (initial condition = $(\pp_{\text{init}}, \qq_{\text{init}})$).}
    \label{fig:effects_of_sequence_len_H0=2.00003}
\end{figure}

\begin{figure}[h!]
    \centering
    \begin{subfigure}[b]{0.48\linewidth}
        \centering
        \includegraphics[width=\linewidth]{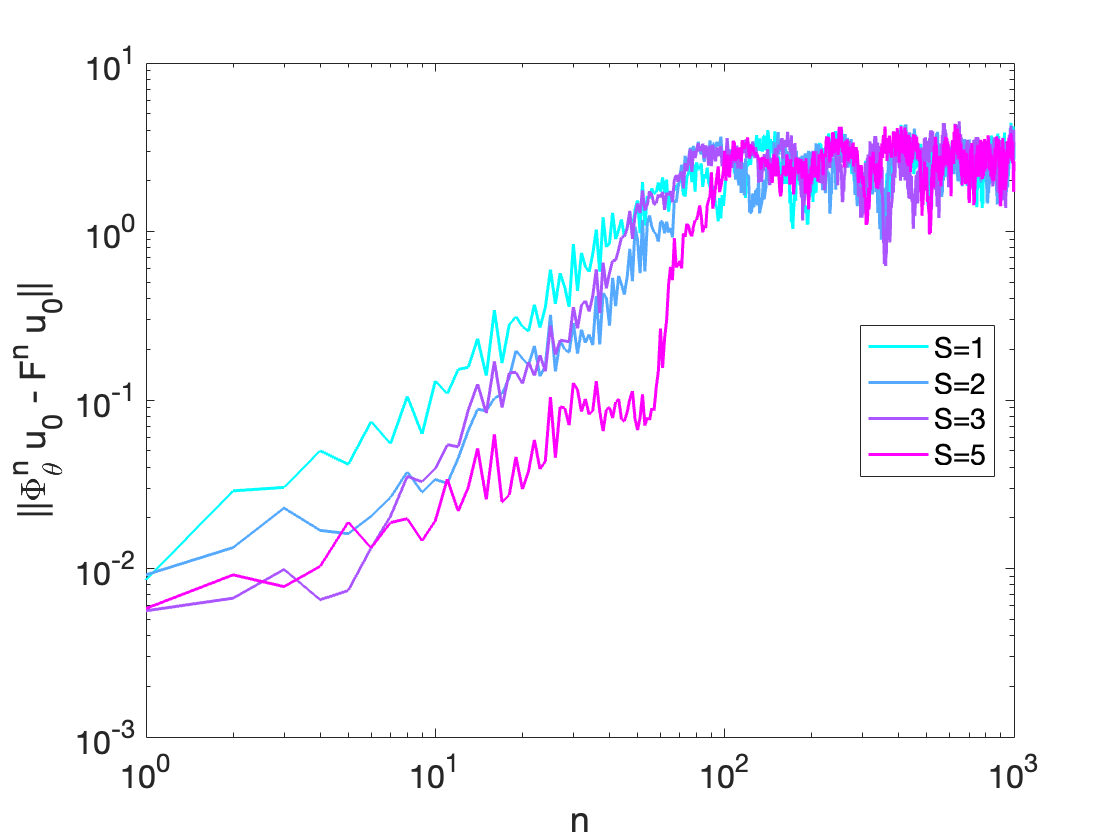}
        \caption{trajectory error}
    \end{subfigure}
    \hfill
    \begin{subfigure}[b]{0.48\linewidth}
        \centering
        \includegraphics[width=\linewidth]{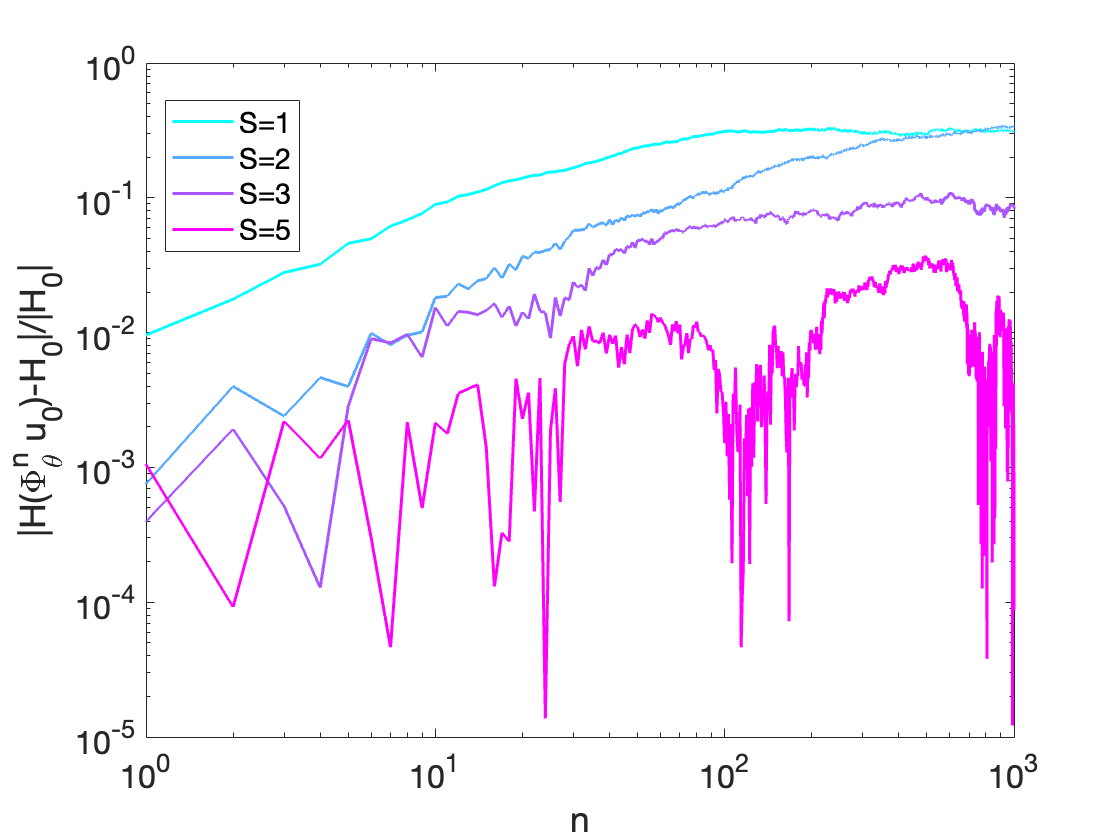}
        \caption{energy error}
    \end{subfigure}
    \caption{Errors in trajectories generated by NN solvers learned with different sequence length $S$ (initial condition = $( \sqrt{2} \pp_{\text{init}}, \qq_{\text{init}})$).}
    \label{fig:effects_of_sequence_len_H0=3.00003}
\end{figure}

We trained the shallow network ResNet(4, 1000) using $\mathcal{D}^{\text{HMC-$H_0$}}$ and multi-step MSE loss for different sequence length $S$. Figure~\ref{fig:effects_of_sequence_len_H0=2.00003} shows that longer sequence length yields slightly better accuracy for the first 10 steps. After 10 steps, there is no significant difference between results of different sequence lengths. In Figure~\ref{fig:effects_of_sequence_len_H0=3.00003}, we repeated the comparison for a different initial condition $(\sqrt{2} \pp_{\text{init}}, \qq_{\text{init}})$. Note that the corresponding energy level is higher than $H_0$ for generating the training data. In other words, we are testing the generalization ability of the NN solvers for out-of-distribution examples. Based on the results, the NN solvers are able to achieve accuracy on par with the accuracy for in-distribution examples for at least the first few steps. Moreover, we found longer training sequences result in significantly better generalization ability.

\subsubsection*{Effects of loss function metric}

\begin{figure}[h!]
    \centering
    \begin{subfigure}[b]{0.48\linewidth}
        \centering
        \includegraphics[width=\linewidth]{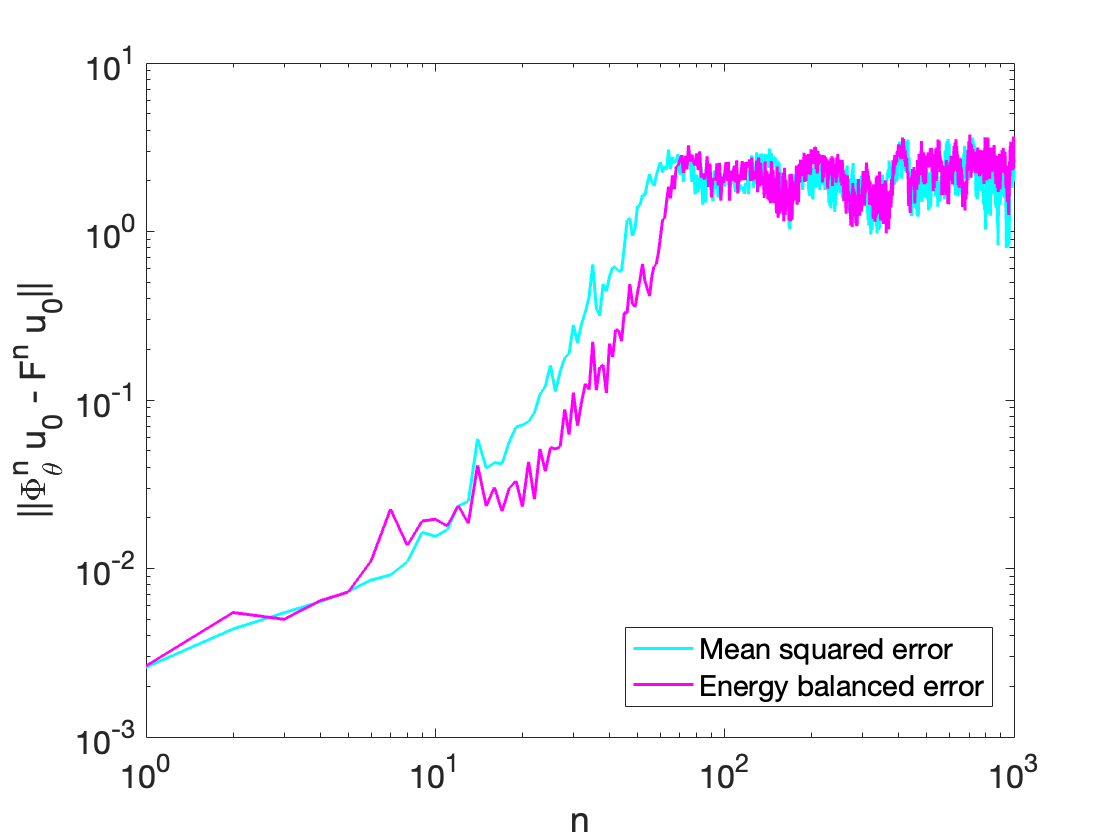}
        \caption{trajectory error}
    \end{subfigure}
    \hfill
    \begin{subfigure}[b]{0.48\linewidth}
        \centering
        \includegraphics[width=\linewidth]{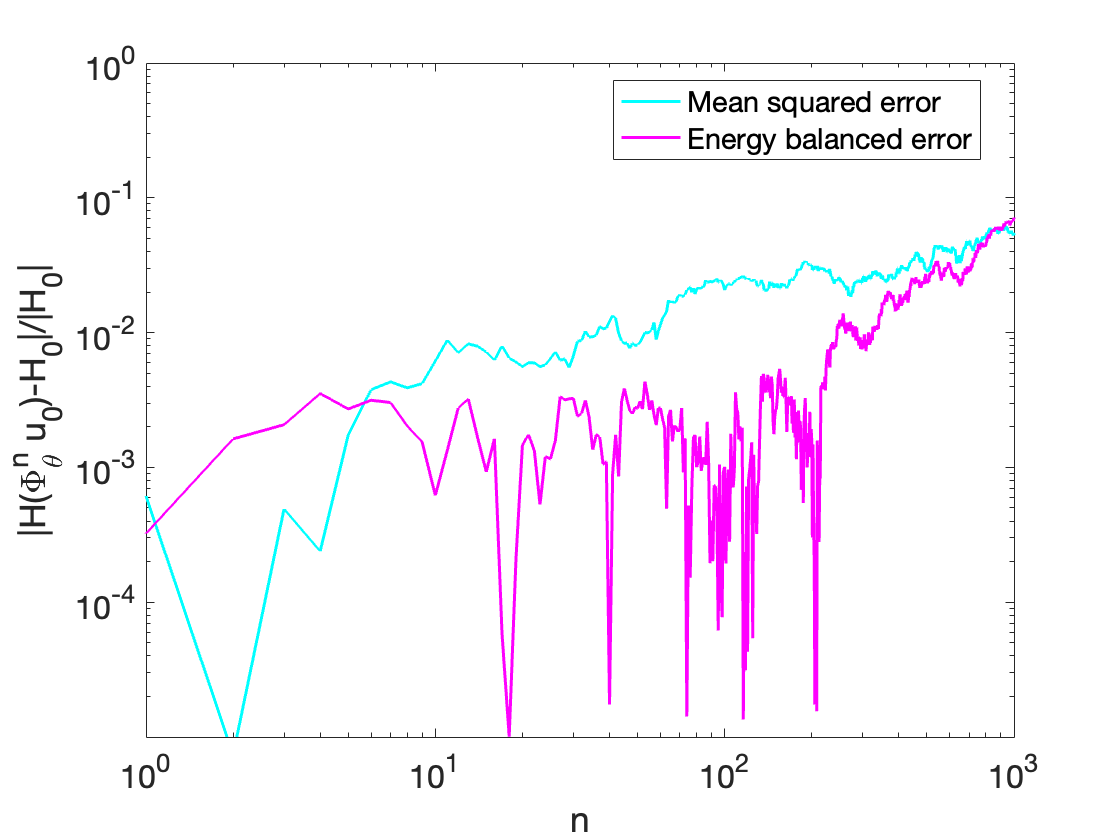}
        \caption{energy error}
    \end{subfigure}
    \caption{Errors in trajectories generated by NN solvers learned with different metrics.}
    \label{fig:effects_of_metric}
\end{figure}

We trained the shallow network ResNet(4, 1000) using $\mathcal{D}^{\text{HMC-$H_0$}}$ and different loss metrics with sequence length $S=5$. As displayed in Figure~\ref{fig:effects_of_metric}, compared to using MSE, using EBE leads to smaller trajectory error and energy error for over 100 steps. In particular, the energy error is not only smaller but also more stable over a long time period.

% \subsubsection*{Effects of network depth}

% % remove this section 
% % say in one sentence .... 
% Lastly, we compare the performance of ResNet(4, 1000) and ResNet(75, 200). We used $\mathcal{D}^{\text{HMC-$H_0$}}$ and the multi-step ($S=5$) EBE loss. Figure~\ref{fig:effects_of_depth} shows ResNet(75, 200) results in lower trajectory error within the first 20 steps, but later the trajectory error increases faster than that of the ResNet(4, 1000) solution. On the other hand, the energy errors of both solutions are small. The energy error of the ResNet(75, 200) solution stays low for a longer period of time, suggesting that deeper networks may produce more stable solutions. We shall continue the comparison in Section 4.4 through evaluating how the NN solution maps can be improved by parareal iterations.

% \begin{figure}[h!]
%     \centering
%     \begin{subfigure}[b]{0.48\linewidth}
%         \centering
%         \includegraphics[width=\linewidth]{figs/nn/depth/traj_err_loglog.png}
%         \caption{trajectory error}
%     \end{subfigure}
%     \hfill
%     \begin{subfigure}[b]{0.48\linewidth}
%         \centering
%         \includegraphics[width=\linewidth]{figs/nn/depth/H_err_loglog.png}
%         \caption{energy error}
%     \end{subfigure}
%     \caption{Errors in the generated trajectories by NN solution maps learned with different network architectures.}
%     \label{fig:effects_of_depth}
% \end{figure}

\subsubsection*{Comparison with numerical solvers}

We present in Table~\ref{tab:fpu-solnmap} the one-step accuracy and runtime performance of various solvers. The fine solver $F_{\Delta t}$ is $\Phi^{\text{KL8}, h=2^{-18}}_{\Delta t}$ implemented in quadruple precision. The NN solver $\Phi^{\text{NN}}_{\Delta t}$ is ResNet(4, 1000), trained using $\mathcal{D}^{\text{HMC-$H_0$}}$ and the multi-step ($S=5$) EBE loss. For comparison, we include several numerical solvers: $\Phi^{\text{CSS4}, h=2^{-9}}_{\Delta t}$ and $\Phi^{\text{VV}, h=2^{-14}}_{\Delta t}$, whose one-step trajectory error is comparable to that of $\Phi^{\text{NN}}_{\Delta t}$, as well as $\Phi^{\text{CSS4}, h=2^{-8}}_{\Delta t}$ and $\Phi^{\text{VV}, h=2^{-11}}_{\Delta t}$, whose one-step energy error is comparable to that of $\Phi^{\text{NN}}_{\Delta t}$. Here VV stands for the \nth{2}-order velocity Verlet scheme. 

It can be seen that, with the same level of trajectory error, the numerical solvers achieve lower energy error than $\Phi^{\text{NN}}_{\Delta t}$. However, in terms of runtime, $\Phi^{\text{NN}}_{\Delta t}$ is as good as $\Phi^{\text{CSS4}, h=2^{-9}}_{\Delta t}$ and is about 17 times faster than $\Phi^{\text{VV}, h=2^{-14}}_{\Delta t}$. We emphasize that the runtime measurements took place in different environments: the NN solver is implemented in Python and numerical solvers are implemented in Julia. We have fully optimized the Julia code for runtime and memory efficiency. We expect to further optimize the NN implementation for better runtime performance in the future.

\begin{table}[h!]
    \caption{Accuracy and runtime performance comparison for various solvers.  *Runtime for the NN solver implemented in Python. The rest of the solvers are written and optimized in Julia.}
    \label{tab:fpu-solnmap}
    \bgroup 
    \def\arraystretch{1.5}
    \begin{tabular}{llll}
        \toprule
        & trajectory error & energy error & runtime \\  
        \midrule 
        $F_{\Delta t}$ & 0 & 0 & 12.3 s \\
        $\Phi^{\text{NN}}_{\Delta t}$ & 0.00264 & $3.2 \times 10^{-4}$  & 0.272 ms*\\
        $\Phi^{\text{CSS4}, h=2^{-9}}_{\Delta t}$ & 0.00262 & $4.3 \times 10^{-6}$ & 0.246 ms \\ 
        $\Phi^{\text{CSS4}, h=2^{-8}}_{\Delta t}$ & 0.0435 & $1.7 \times 10^{-4}$ & 0.150 ms \\ 
        $\Phi^{\text{VV}, h=2^{-14}}_{\Delta t}$ & 0.00419 & $3.4 \times 10^{-7}$ & 4.231 ms \\ 
        $\Phi^{\text{VV}, h=2^{-11}}_{\Delta t}$ & 0.231 & $6.5 \times 10^{-4}$ & 0.578 ms \\ 
        \botrule
    \end{tabular}
    \egroup 
\end{table}

\subsection{NN solution map in parareal iterations}

In this section, we present results of using $\Phi^{\text{NN}}_{\Delta t}$ as the coarse solver in parareal methods.

\begin{figure}[h!]
    \centering
    \begin{tabular}{ccc}
        &  Traj error  & Energy error  \\
        $\Phi^{\text{NN}}_{\Delta t}$ & \adjincludegraphics[height=4cm, valign=m, trim={{0.04\width} 0 {0.06\width} 0}, clip]{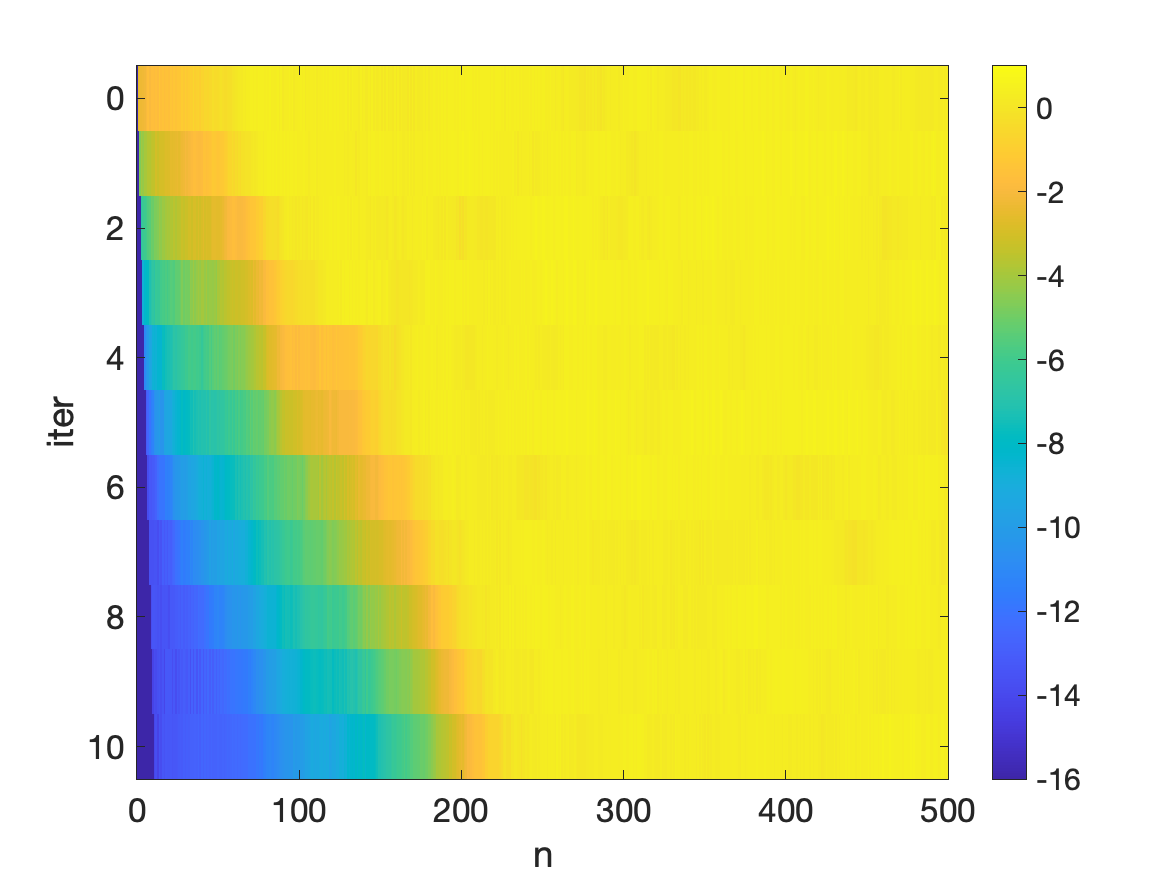} & \adjincludegraphics[height=4cm, valign=m, trim={{0.04\width} 0 {0.06\width} 0}, clip]{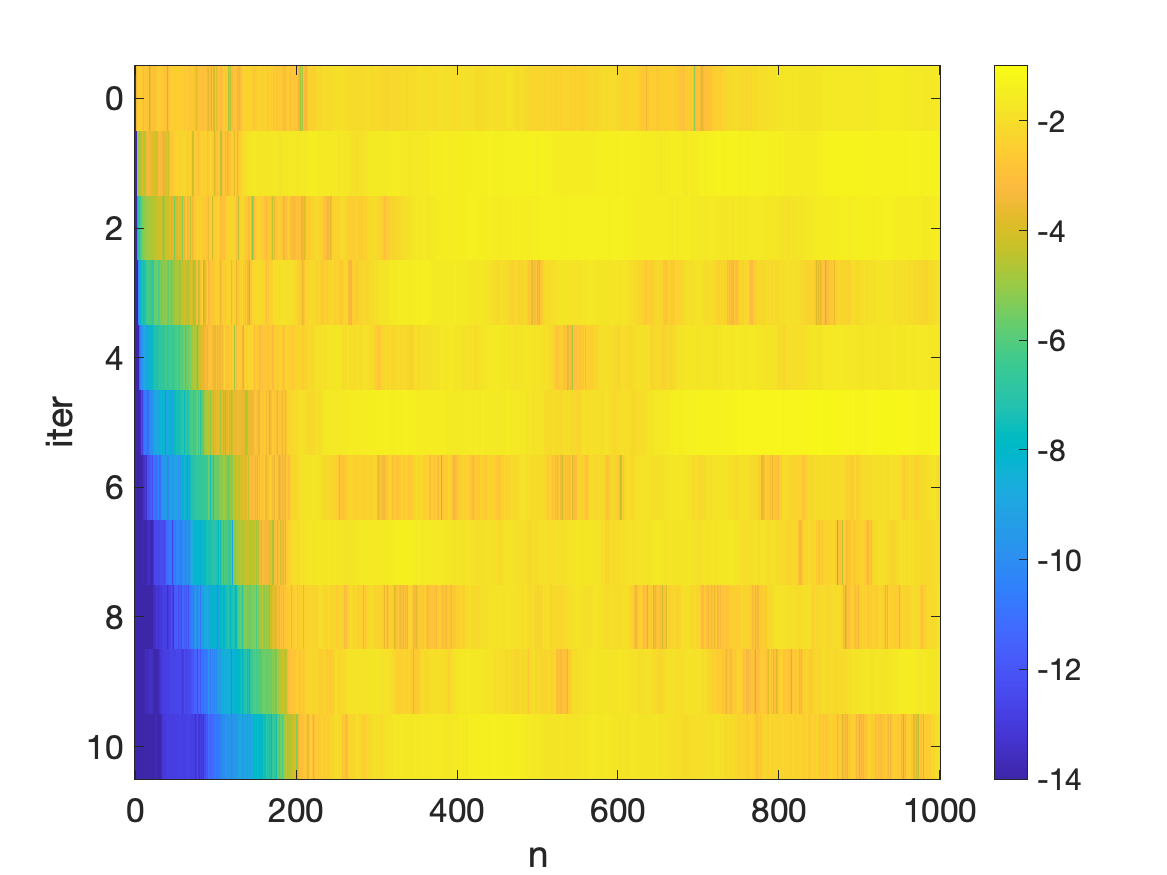} \\
        $\Phi^{\text{CSS4}, h=2^{-9}}_{\Delta t}$ &  \adjincludegraphics[height=4cm, valign=m, trim={{0.04\width} 0 {0.06\width} 0}, clip]{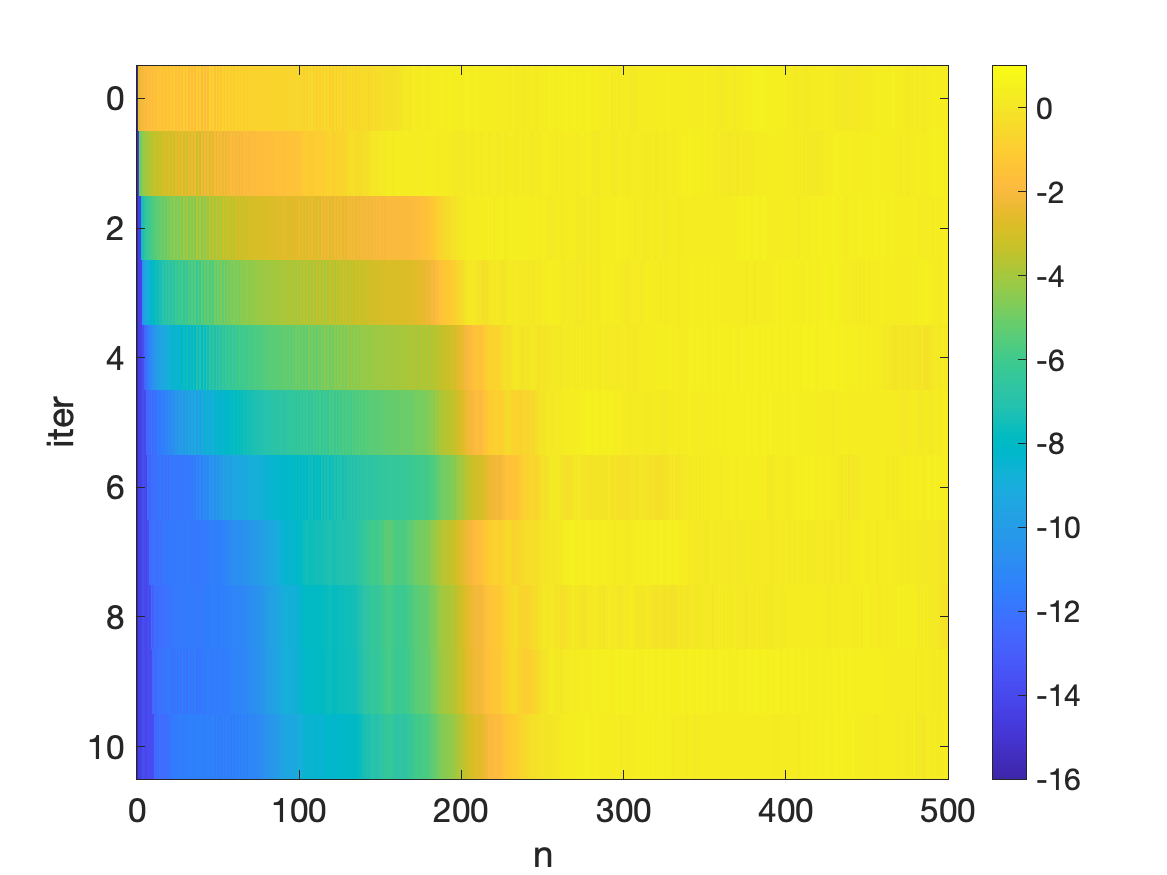} & \adjincludegraphics[height=4cm, valign=m, trim={{0.04\width} 0 {0.06\width} 0}, clip]{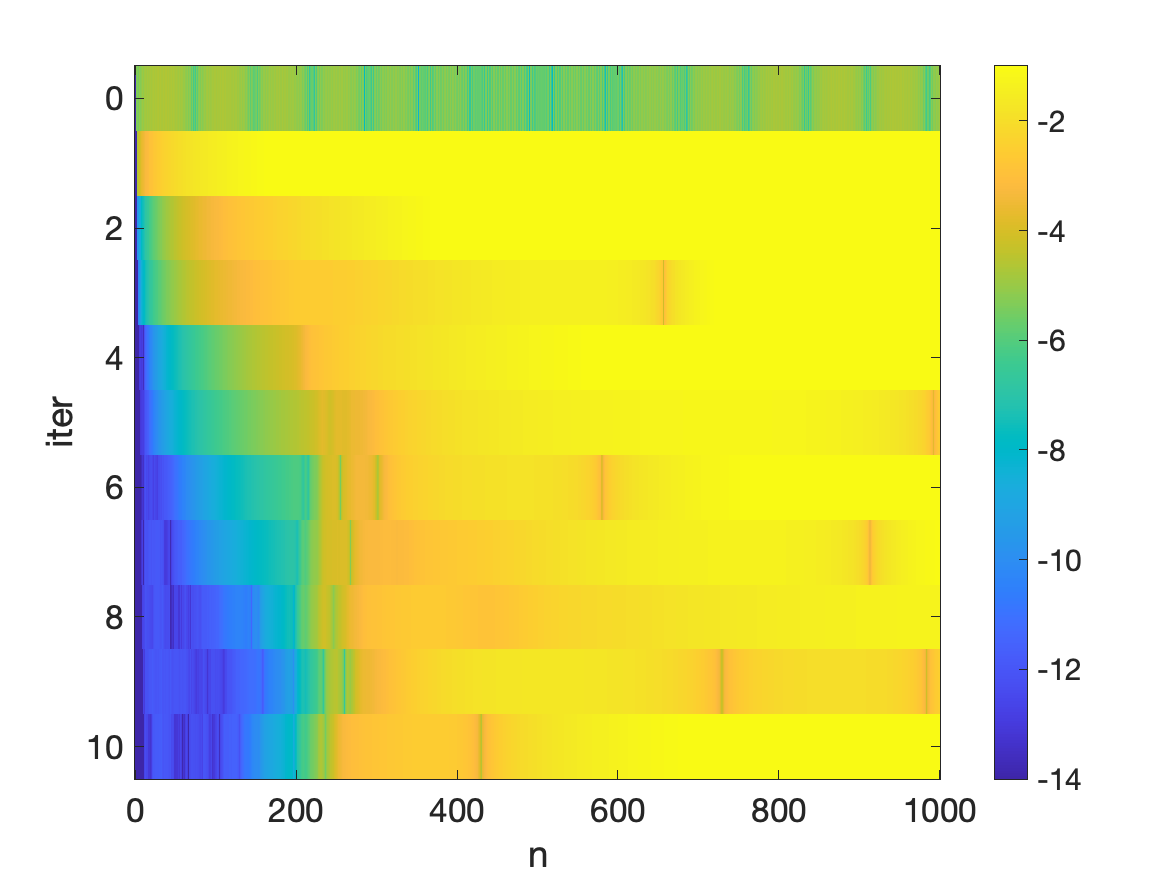} \\
        $\Phi^{\text{CSS4}, h=2^{-8}}_{\Delta t}$ & \adjincludegraphics[height=4cm, valign=m, trim={{0.04\width} 0 {0.06\width} 0}, clip]{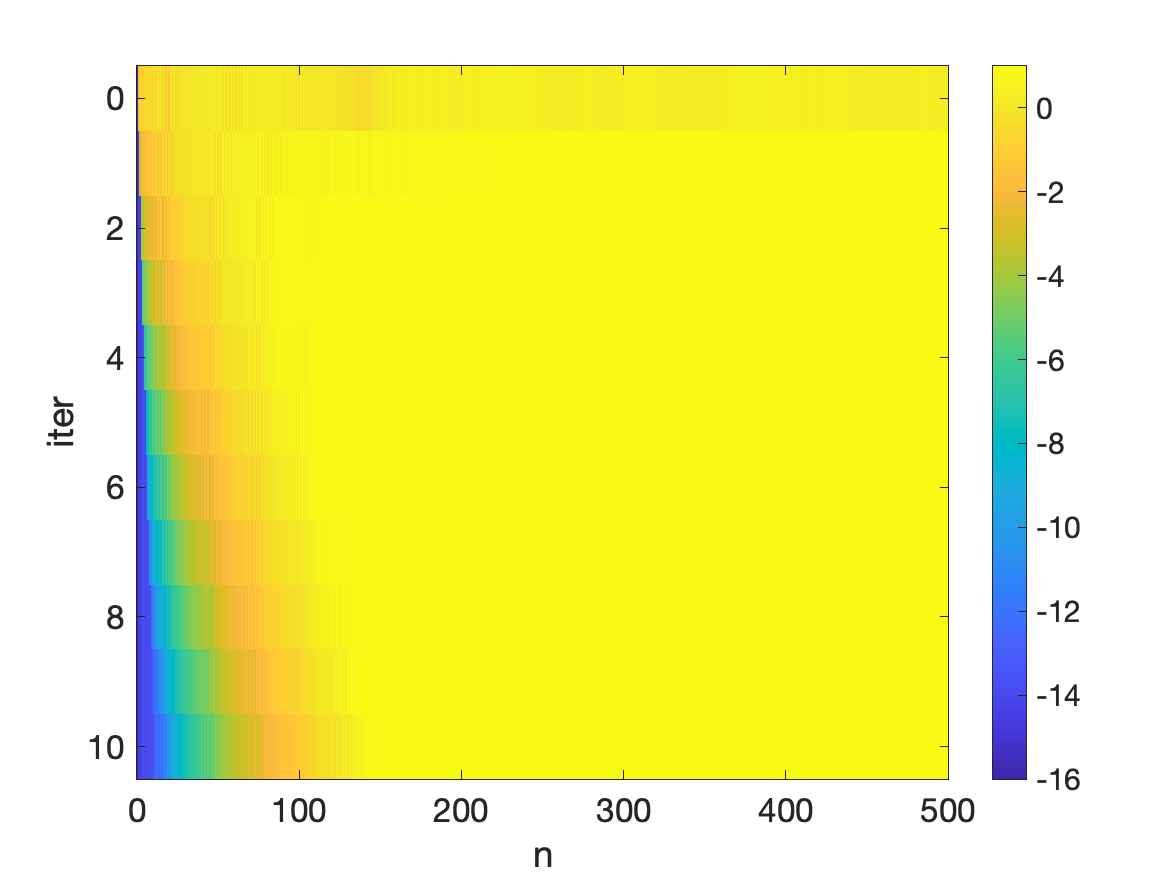}  & \adjincludegraphics[height=4cm, valign=m, trim={{0.04\width} 0 {0.06\width} 0}, clip]{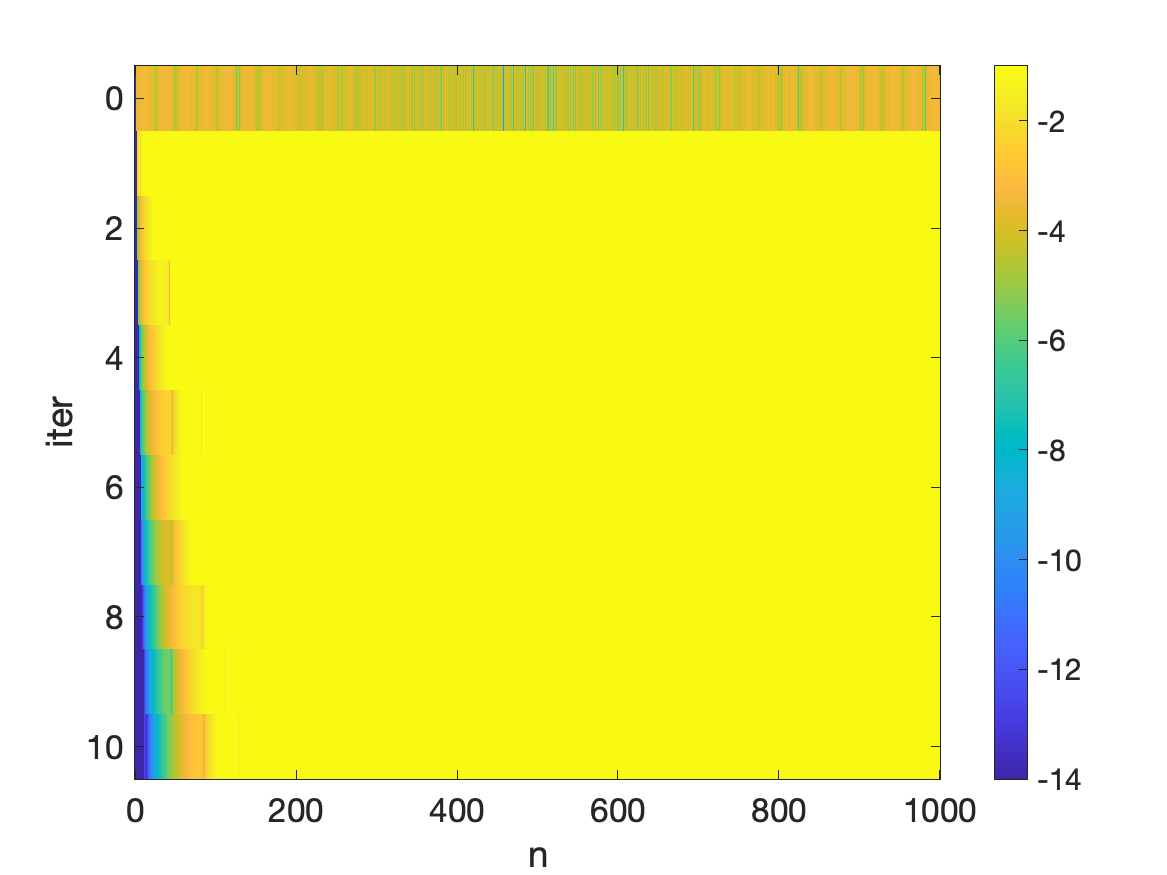} 
    \end{tabular}
    \caption{Log (base 10) errors in plain parareal  solutions by various coarse solvers ($\Delta t=1.0$, $T=1000$).}
    \label{fig:plain-parareal-nn}
\end{figure}

\begin{figure}[h!]
    \centering
    \begin{tabular}{ccc}
        &  $k=0$ & $k=3$ \\
        $\Phi^{\text{NN}}_{\Delta t}$ & 
        \adjincludegraphics[height=4cm, valign=m, trim={{0.04\width} 0 {0.06\width} 0}, clip]{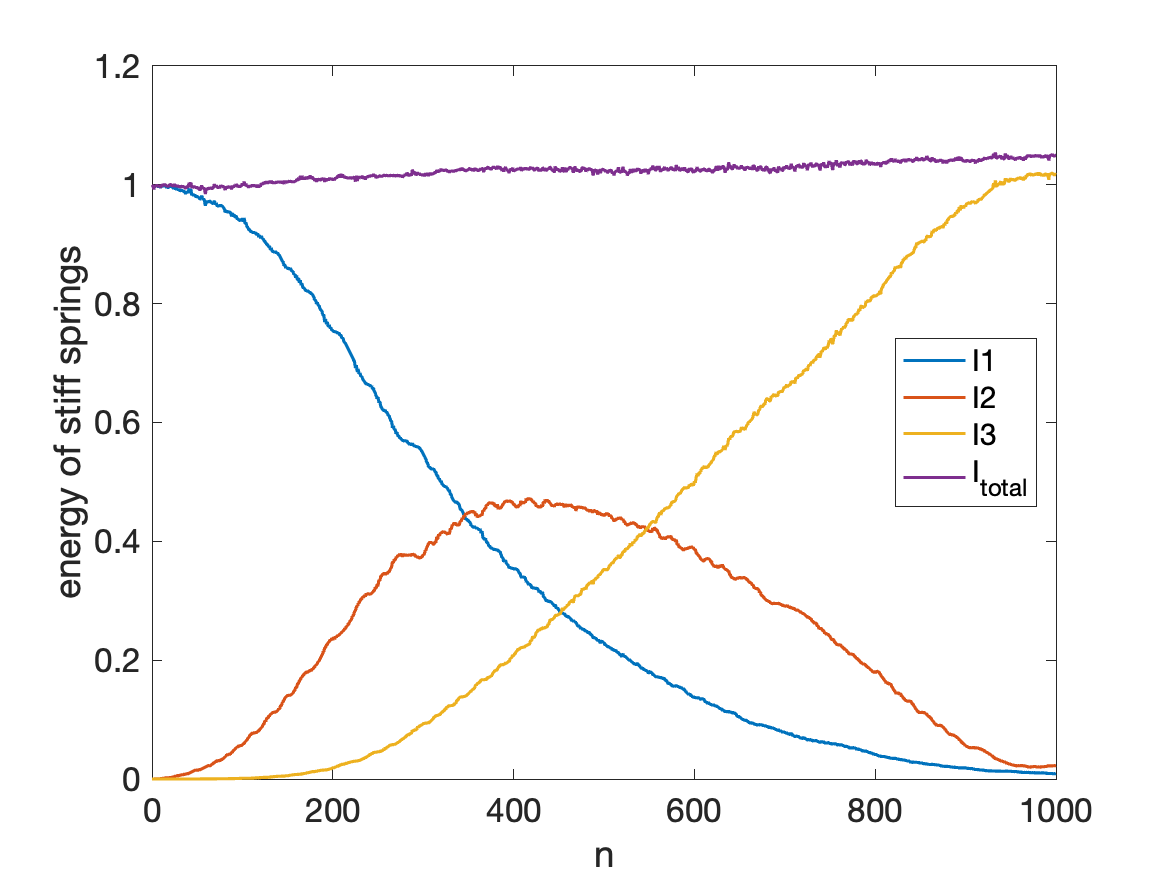} & \adjincludegraphics[height=4cm, valign=m, trim={{0.04\width} 0 {0.06\width} 0}, clip]{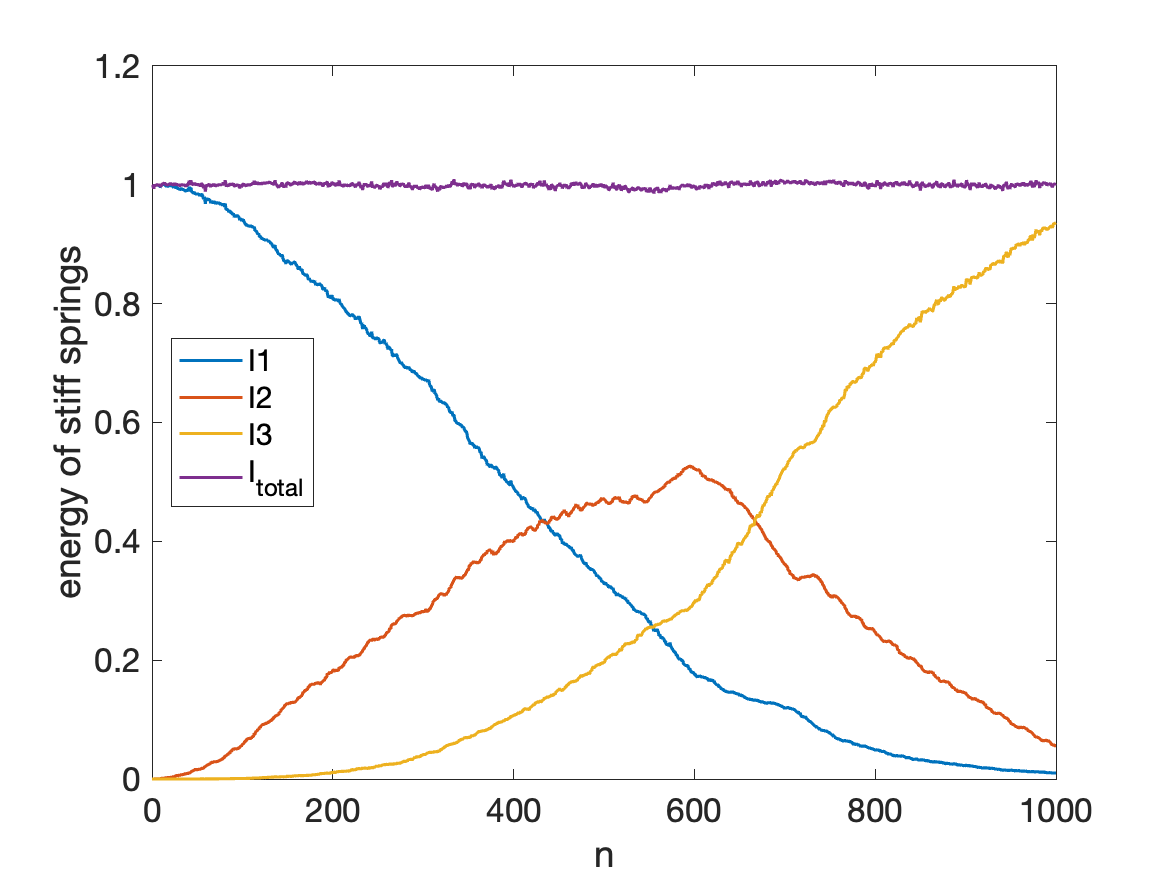} \\
         $\Phi^{\text{CSS4}, h=2^{-9}}_{\Delta t}$ & \adjincludegraphics[height=4cm, valign=m, trim={{0.09\width} 0 {0.06\width} 0}, clip]{figs/parareal/plain/css4-h=2e-9/k=0-stiffenergies.png} & 
        \adjincludegraphics[height=4cm, valign=m, trim={{0.09\width} 0 {0.06\width} 0}, clip]{figs/parareal/plain/css4-h=2e-9/k=3-stiffenergies.png} 
        \\
        $\Phi^{\text{CSS4}, h=2^{-8}}_{\Delta t}$ & \adjincludegraphics[height=4cm, valign=m, trim={{0.09\width} 0 {0.06\width} 0}, clip]{figs/parareal/plain/css4-h=2e-8/k=0-stiffenergies.png} & \adjincludegraphics[height=4cm, valign=m, trim={{0.09\width} 0 {0.06\width} 0}, clip]{figs/parareal/plain/css4-h=2e-8/k=3-stiffenergies.png} 
    \end{tabular}
    \caption{Energy profiles of the stiff springs computed from plain parareal solutions by various coarse solvers ($\Delta t=1.0$, $T=1000$).}
    \label{fig:plain-parareal-nn-stiffenergies}
\end{figure}

We first study the plain parareal method. Figure~\ref{fig:plain-parareal-nn} compares the plain parareal solutions computed by different coarse solvers, including $\Phi^{\text{NN}}_{\Delta t}$, $\Phi^{\text{CSS4}, h=2^{-9}}_{\Delta t}$ and $\Phi^{\text{CSS4}, h=2^{-8}}_{\Delta t}$. Clearly, $\Phi^{\text{CSS4}, h=2^{-8}}_{\Delta t}$ performs the worst, as expected because it is the least accurate among the three solvers. Based on the trajectory errors, we see using $\Phi^{\text{NN}}_{\Delta t}$ as the coarse solver provides slower accuracy improvement over iterations compared to using $\Phi^{\text{CSS4}, h=2^{-9}}_{\Delta t}$. Comparing the energy errors, we observe that when using $\Phi^{\text{NN}}_{\Delta t}$, the stability in energy is not destroyed as much as in using $\Phi^{\text{CSS4}, h=2^{-9}}_{\Delta t}$ as the coarse solver (see also the energy profiles at iteration~3 in Figure~\ref{fig:plain-parareal-nn-stiffenergies}). 
% This shows that, through replacing the numerical solver by an NN solver of similar accuracy, we are able to stabilize the resulting parareal solutions.

% [how so?]

\begin{figure}[h!]
    \centering
    \begin{tabular}{ccc}
        &  Traj error  & Energy error  \\
        $\Phi^{\text{NN}}_{\Delta t}$ & \adjincludegraphics[height=4cm, valign=m, trim={{0.04\width} 0 {0.06\width} 0}, clip]{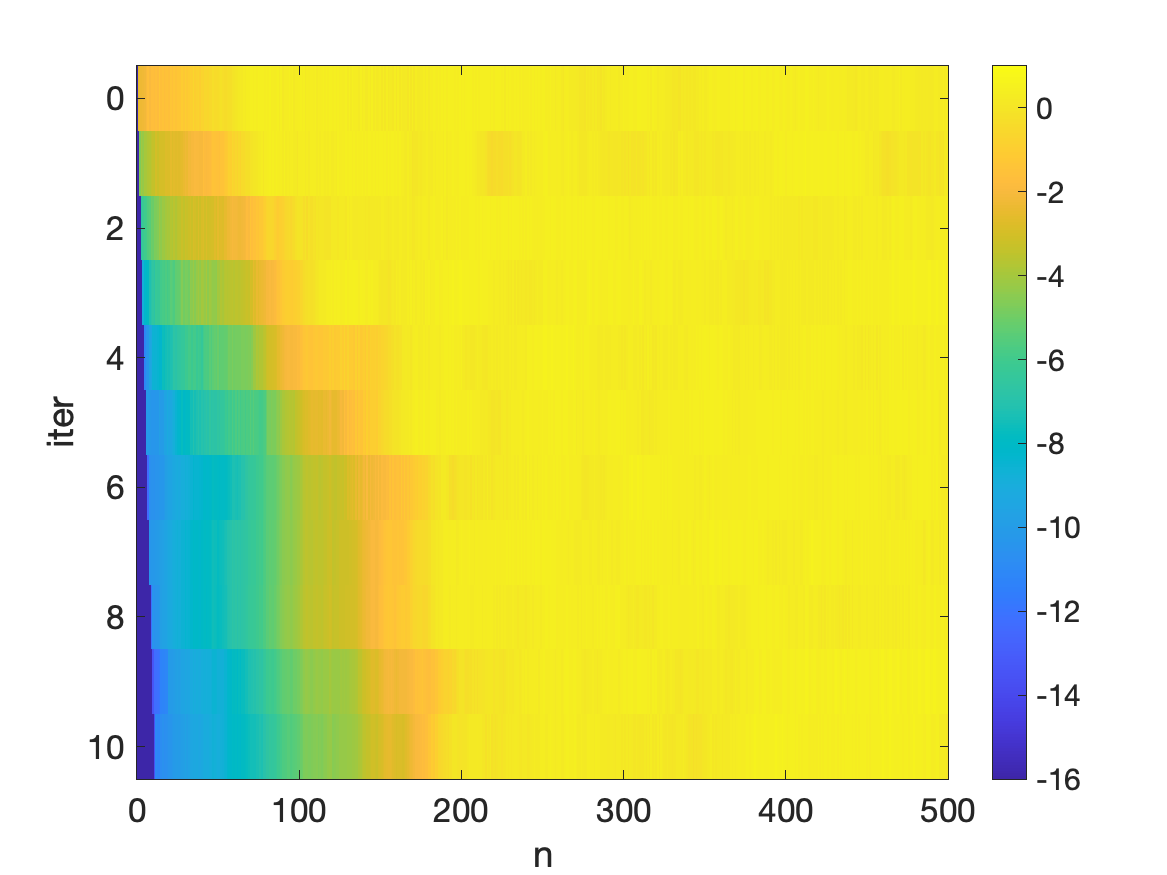} & \adjincludegraphics[height=4cm, valign=m, trim={{0.04\width} 0 {0.06\width} 0}, clip]{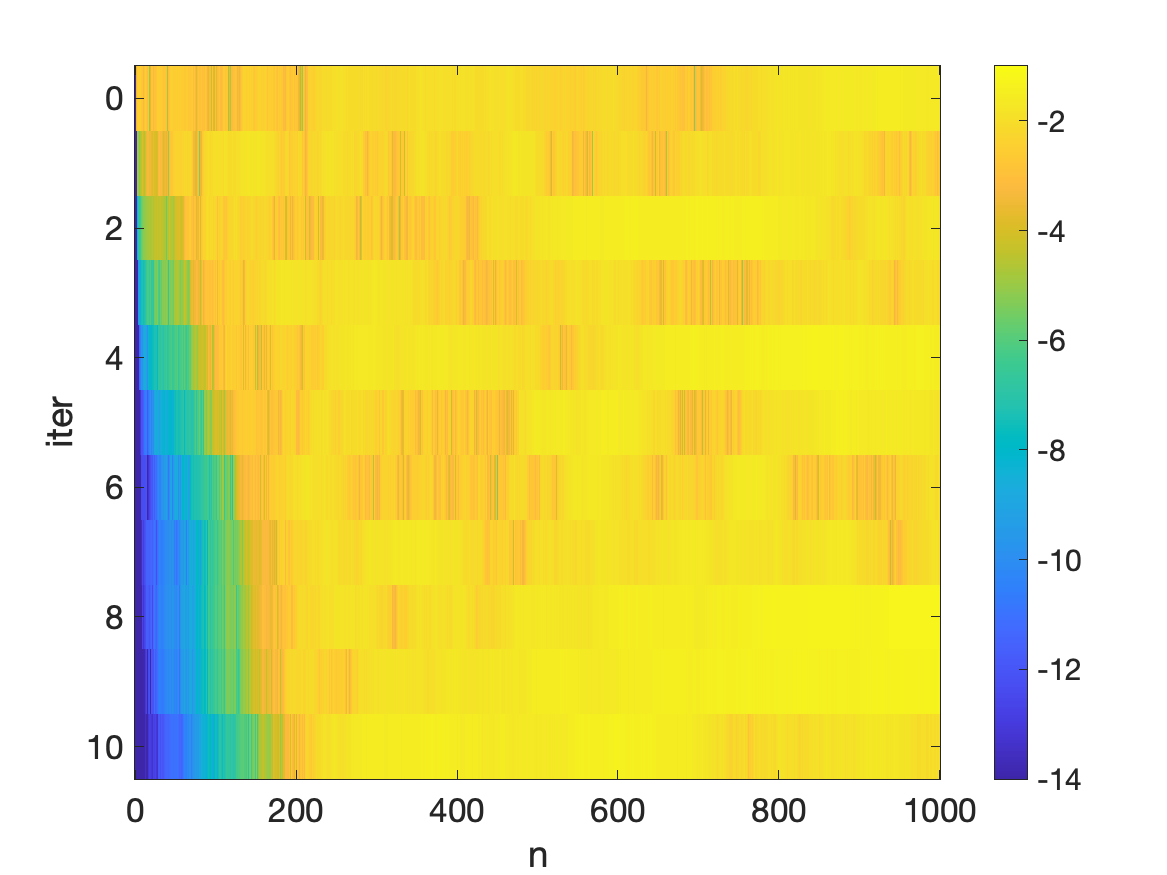} \\
        $\Phi^{\text{CSS4}, h=2^{-9}}_{\Delta t}$ &  \adjincludegraphics[height=4cm, valign=m, trim={{0.04\width} 0 {0.06\width} 0}, clip]{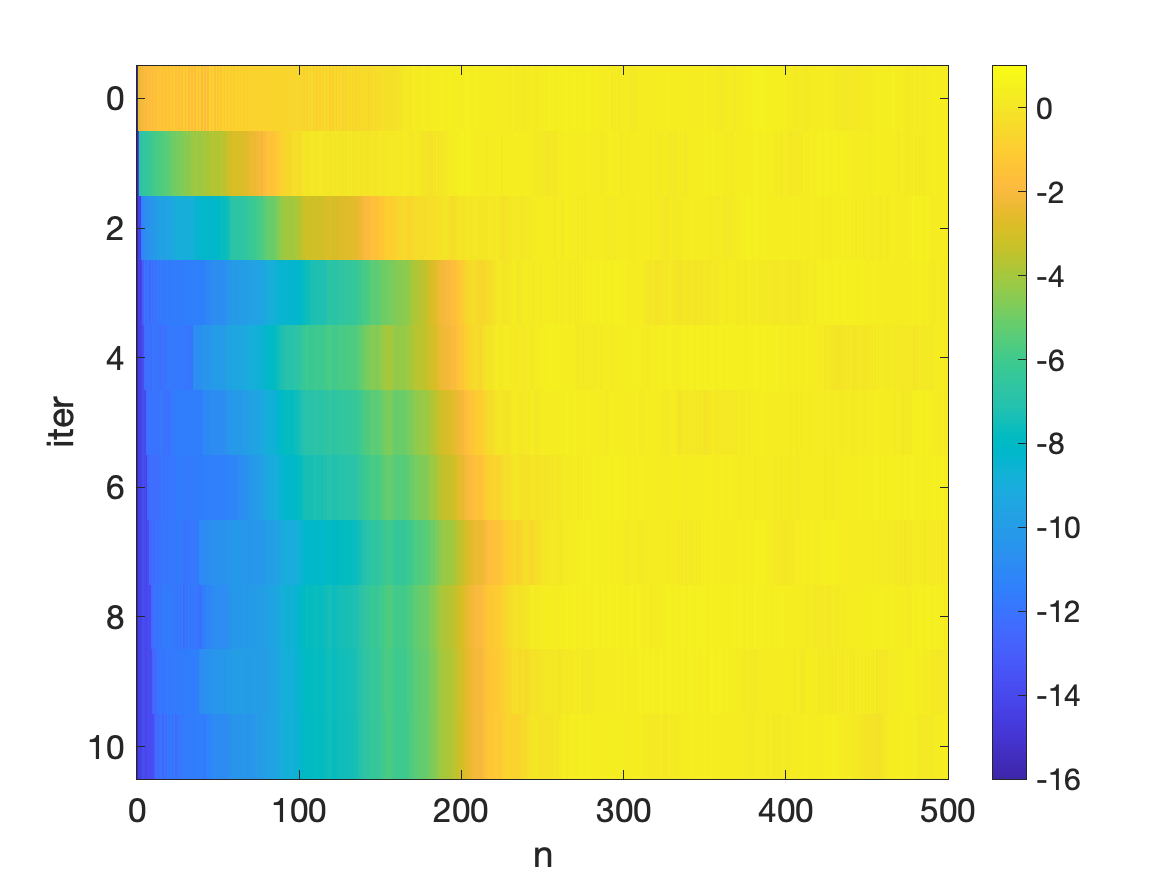} & \adjincludegraphics[height=4cm, valign=m, trim={{0.04\width} 0 {0.06\width} 0}, clip]{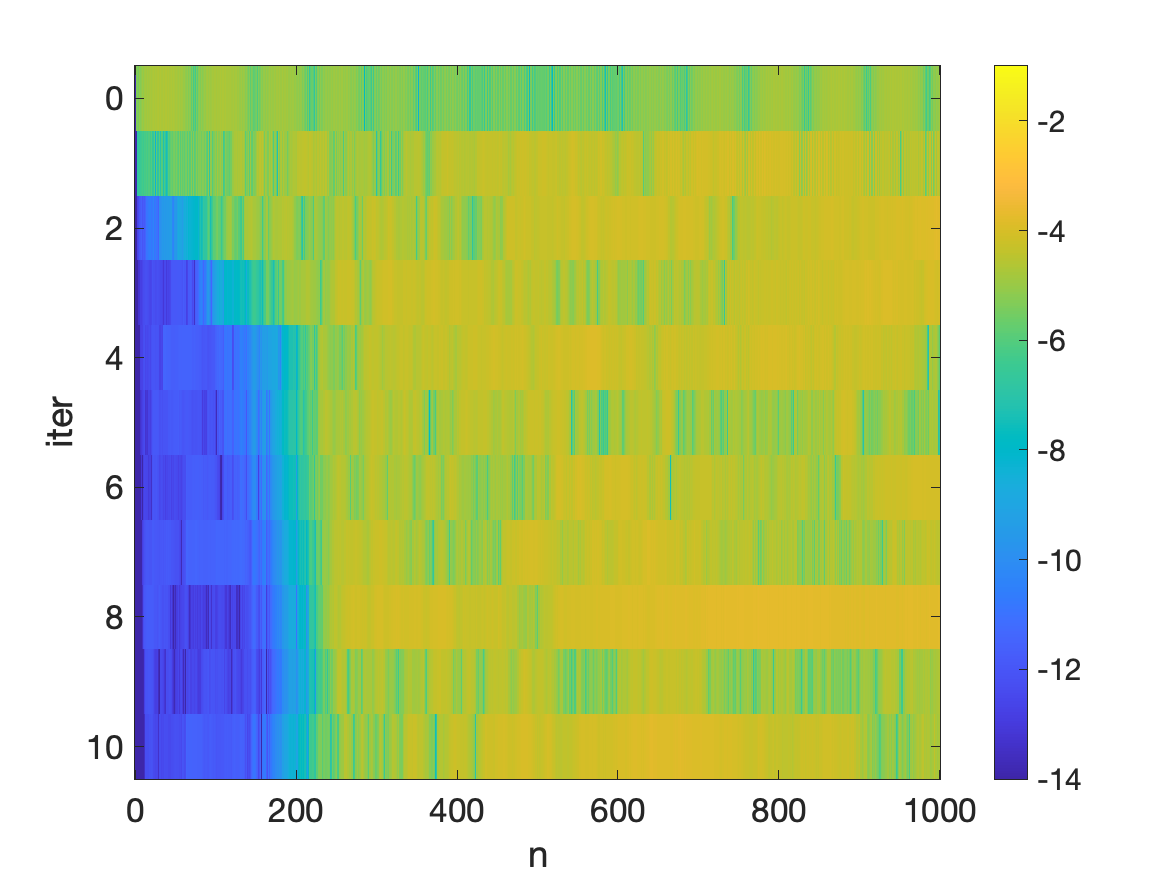} \\
        $\Phi^{\text{CSS4}, h=2^{-8}}_{\Delta t}$ & \adjincludegraphics[height=4cm, valign=m, trim={{0.04\width} 0 {0.06\width} 0}, clip]{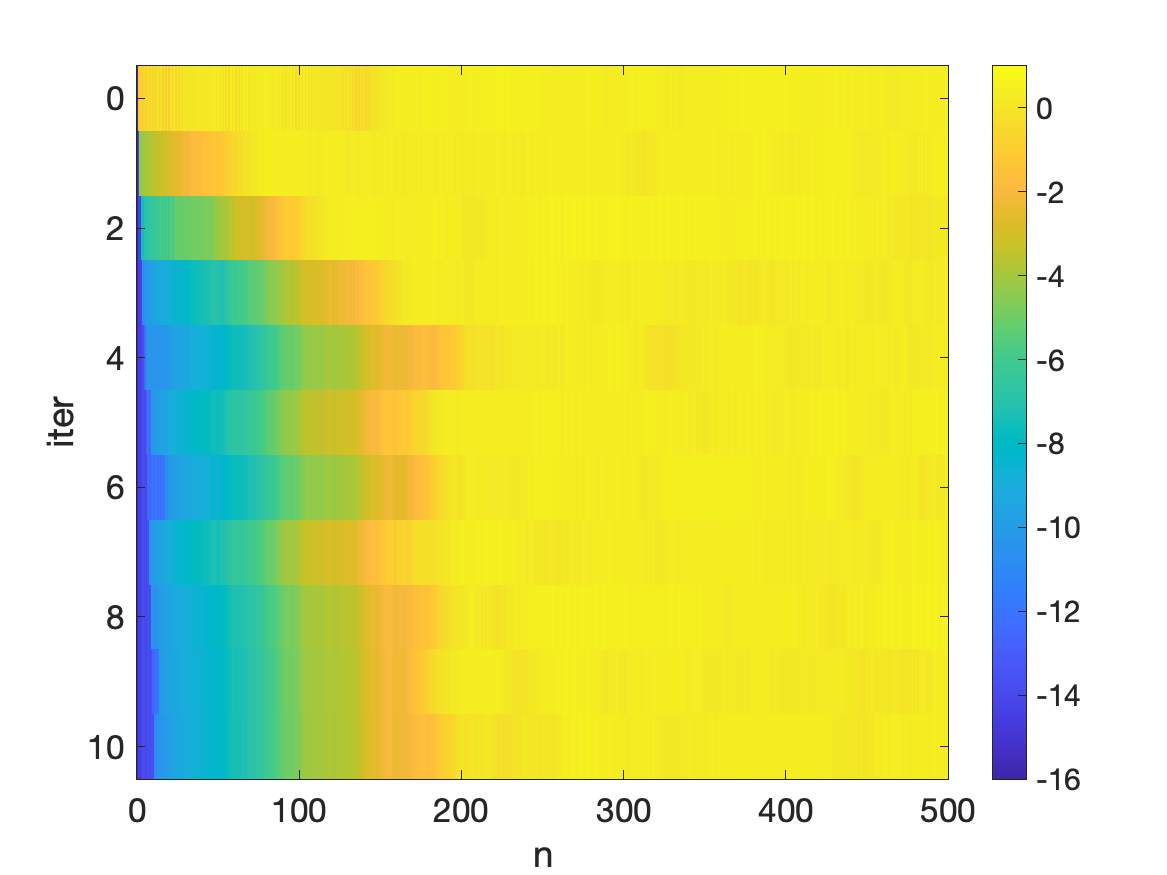}  & \adjincludegraphics[height=4cm, valign=m, trim={{0.04\width} 0 {0.06\width} 0}, clip]{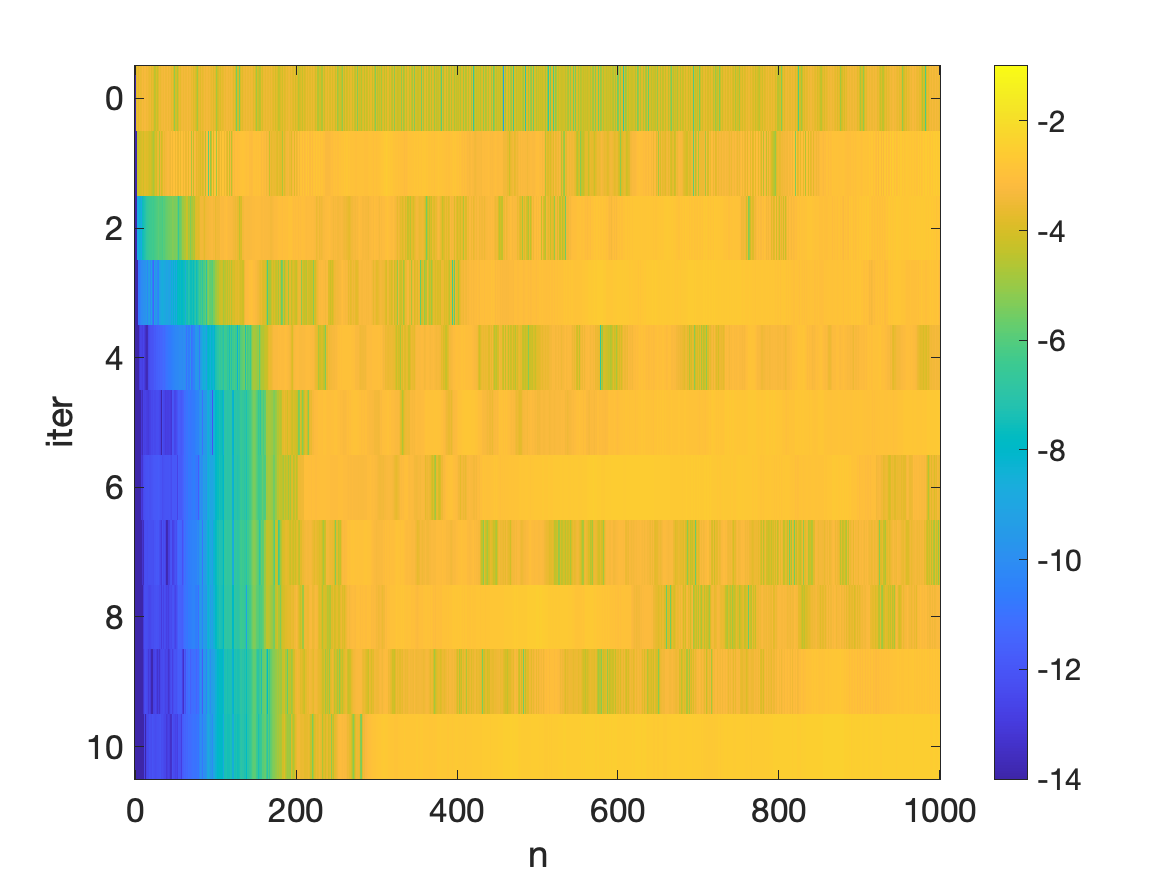} 
    \end{tabular}
    \caption{Log (base 10) errors in Procrustes parareal solutions by various coarse solvers ($\Delta t=1.0$, $T=1000$).}
    \label{fig:procrustes-parareal-nn}
\end{figure}

We will now compare different coarse solvers used in the Procrustes parareal method. In Section 4.2, we found the Procrustes parareal method improves accuracy by stablilizing the energy of the solutions. 
% We would like to study if combining the Procrustes parareal and the NN solution map will be more beneficial. 
As shown in Figure~\ref{fig:procrustes-parareal-nn}, the best improvement is again yielded by $\Phi^{\text{CSS4}, h=2^{-9}}_{\Delta t}$. There is no significant gain from combining Procrustes parareal with the NN solution map. In fact, the improvement over the Procrustes parareal iterations even deteriorates compared to the improvement over the plain parareal iterations.

We speculate that $\Phi^{\text{NN}}_{\Delta t}$ does not perform well in parareal iterations because it was not trained using suitable data. As described in Section 3.3, we sampled training data points from the Liouville density, mainly because it is a natural distribution for learning a solution map to be used in a sequential algorithm. In parareal schemes, since the coarse solver is applied differently than in a sequential algorithm, we would need a different data distribution.  A similar issue has been investigated in \cite{nguyen2023numerical} for approximating the correction operator using NN approach for wave equations. There the authors demonstrated the importance of using training data closer to the ones encountered in the simulations. 

% what is phase error distribution between NN solution and fine solution? 
% does it make sense to correct it 

\begin{figure}[h!]
    \centering
    \begin{tabular}{ccc}
        &  $k=0$ & $k=3$ \\
        $\Phi^{\text{NN}}_{\Delta t}$ & 
        \adjincludegraphics[height=4cm, valign=m, trim={{0.04\width} 0 {0.06\width} 0}, clip]{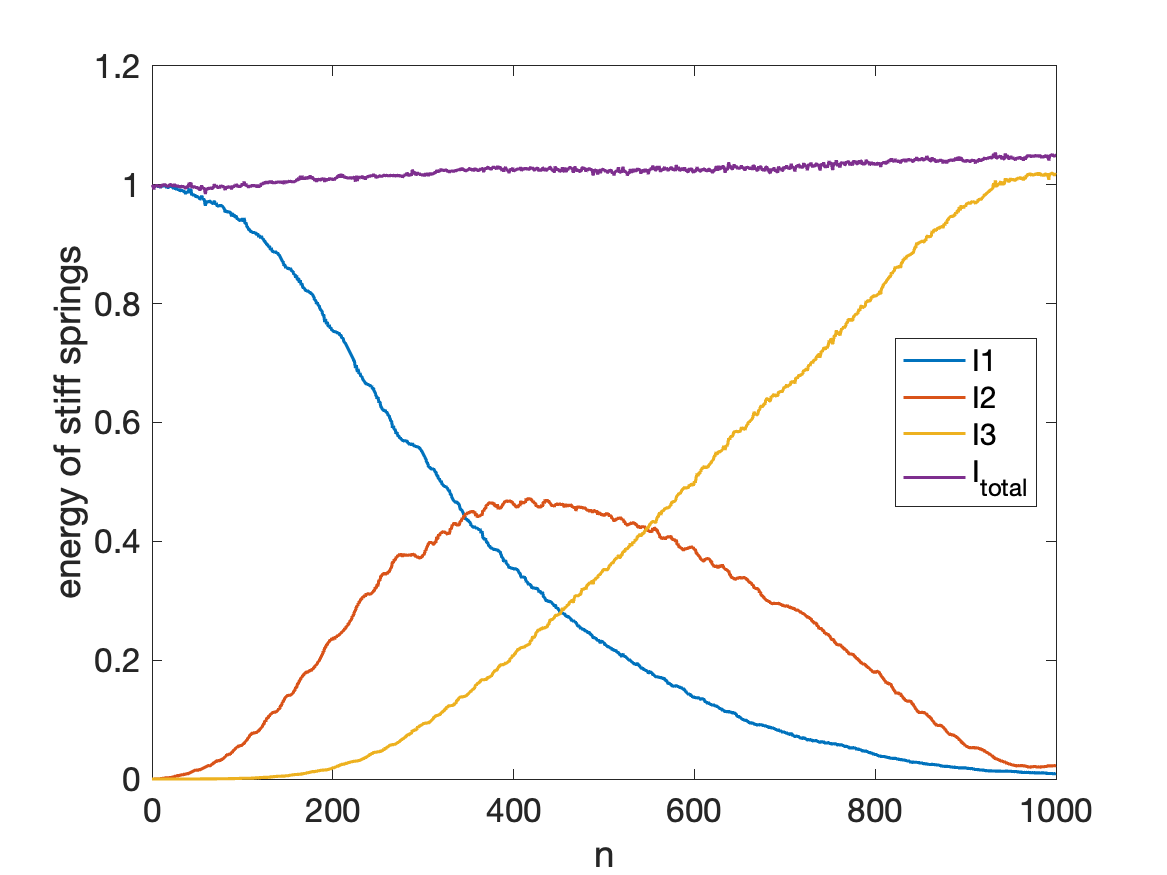} & \adjincludegraphics[height=4cm, valign=m, trim={{0.04\width} 0 {0.06\width} 0}, clip]{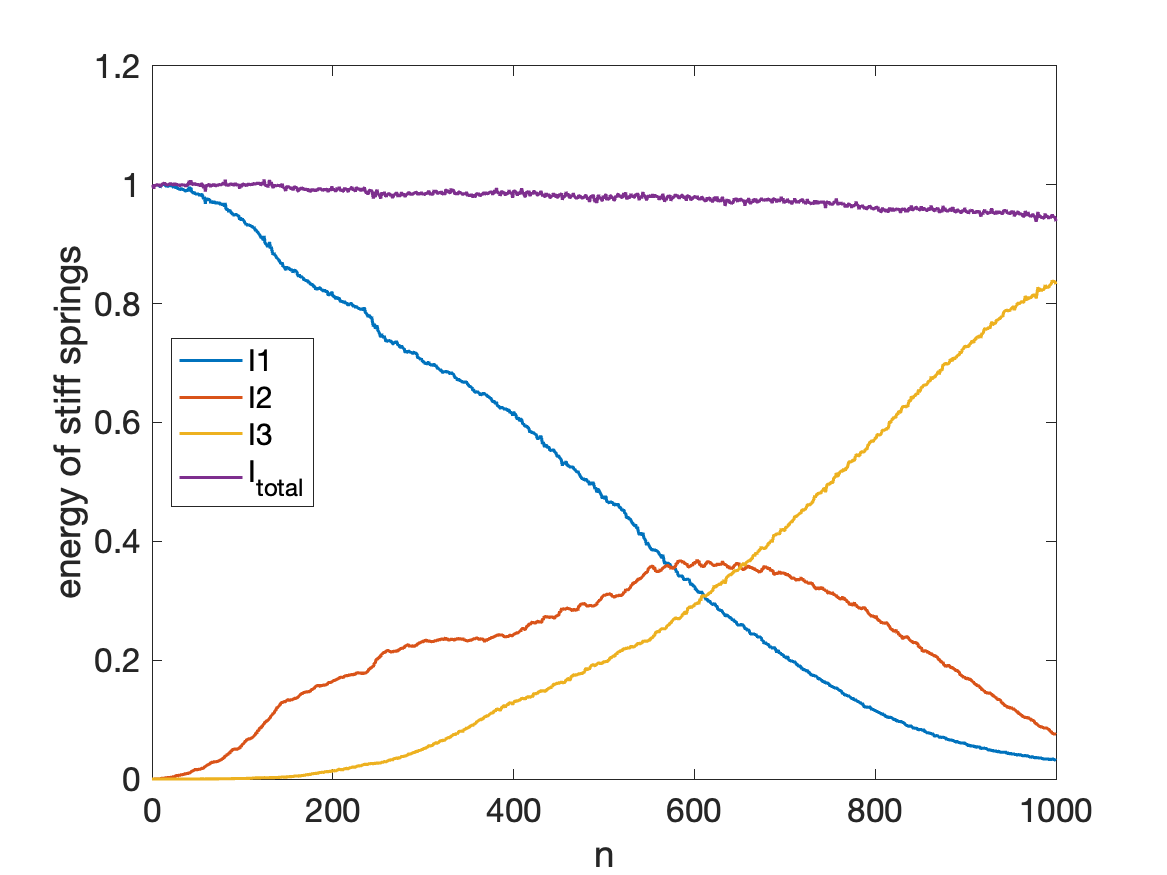} \\
         $\Phi^{\text{CSS4}, h=2^{-9}}_{\Delta t}$ & \adjincludegraphics[height=4cm, valign=m, trim={{0.09\width} 0 {0.06\width} 0}, clip]{figs/parareal/phasecorr/css4-h=2e-9/k=0-stiffenergies.png} & 
        \adjincludegraphics[height=4cm, valign=m, trim={{0.09\width} 0 {0.06\width} 0}, clip]{figs/parareal/phasecorr/css4-h=2e-9/k=3-stiffenergies.png} 
        \\
        $\Phi^{\text{CSS4}, h=2^{-8}}_{\Delta t}$ & \adjincludegraphics[height=4cm, valign=m, trim={{0.09\width} 0 {0.06\width} 0}, clip]{figs/parareal/phasecorr/css4-h=2e-8/k=0-stiffenergies.png} & \adjincludegraphics[height=4cm, valign=m, trim={{0.09\width} 0 {0.06\width} 0}, clip]{figs/parareal/phasecorr/css4-h=2e-8/k=3-stiffenergies.png} 
    \end{tabular}
    \caption{Energy profiles of the stiff springs computed from Procrustes parareal solutions by various coarse solvers ($\Delta t=1.0$, $T=1000$).}
    \label{fig:procrustes-parareal-nn-stiffenergies}
\end{figure}

\section{Conclusion}

In this paper, we presented two data-driven approaches for stabilization of the standard parareal algorithm for long time computation of highly oscillatory Hamiltonian systems. The Procrustes parareal approach uses solutions computed along the parareal iterations to construct a correction operator to align the ``phase'' of the fine and coarse solvers. Numerical results for the FPU problem demonstrated that the constructed correction can successfully stabilize the parareal iterations, which helps improve the accuracy of the computed solutions. The second approach we proposed is to use a neural network (NN) to approximate the reference solution map that advances the given state forward in time by a large fixed time step. We developed a sampling algorithm called HMC-$H_0$ to sample phase space points from the neighborhood of an energy level set. We also designed a loss function which considers the energy balanced errors between approximated trajectories and reference trajectories. The resulting NN solver for the FPU problem is able to achieve comparable or better runtime performance compared to numerical solvers of similar accuracy. When combined with parareal iterations, solutions computed by the NN solver are not as accuracte as solutions computed by a comparable numerical solver, although the NN energy errors are slightly smaller. We think that this may be improved if we train the network using data suitably sampled for the discrete trajectories computed by the parareal schemes. 

The FPU problem is too small to reveal the potential benefit of using NNs. It is small enough that optimized high order symplectic integrators are extremely efficient and accurate. For more complicated problems, where the phase space is large and lower-order accurate methods are the only feasible choice, we think that the investigated NN approach may become viable. 
%One such example is the second order wave equation with discontinuous wave speed as reported in \cite{nguyen2023numerical}.

\backmatter

\section*{Declarations}

\begin{itemize}
\item Ethical approval

Not applicable 

\item Availability of supporting data

The datasets generated and analysed during the current study are available from the corresponding author on reasonable request.

\item Competing interests 

The authors have no competing interests to declare that are relevant to the content of this article.

\item Funding

The authors are partially supported by National Science Foundation grant DMS-2208504.

\item Authors' contributions

R.F.: Conceptualization, Methodology, Software, Visualization, Writing

R.T.: Conceptualization, Funding acquisition, Methodology, Project administration, Supervision, Writing 

\item Acknowledgement

The authors thank the Texas Advanced Computing Center (TACC) for providing computing resources.

\end{itemize}

\noindent

\bibliography{ref}

\end{document}